\documentstyle[12pt]{article}
\begin{document}
\setlength{\textheight}{9in}

\author{S.V. Ludkovsky.}

\title{Skew idempotent functionals of ordered semirings.}

\date{21 December 2012}
\maketitle

\begin{abstract}
Skew idempotent functionals of ordered semirings are studied.
Different associative and non-associative semirings are considered.
Theorems about properties of skew idempotent functionals are proved.
Examples are given. \footnote{  
\par key words: idempotent functional, skew, semiring, quasiring}
\end{abstract}

\section{Introduction.}
Idempotent mathematics being a new branch has attracted much
attention in recent times as a theoretical tool having important
applications in mathematics and quantum physics (see
\cite{litmasshmz69,radcmcu98,zarichizv10} and references therein).
On  the other hand, functionals and measures are also used for
studies of representations of groups and algebras (see, for example,
\cite{nai,fell,hew,lujms147:3:08,lujms150:4:08} and references
therein). Idempotent mathematics arise naturally from the
consideration of the quantization and the Plank constant in physics.
Earlier idempotent functionals associative and commutative relative
to the operation $\bigodot $ and $+$ were investigated and
particularly on spaces of continuous real-valued functions on
compact Hausdorff spaces (see
\cite{litmasshmz69,radcmcu98,zarichizv10} and references therein).
\par Apart from previous works, in this article skew idempotent
functionals generally non-associative and noncommutative relative to
the operation $\bigodot $ are investigated (see Section 3). The
axiomatic in some respect is different from the real case. Moreover,
the compactness condition of topological spaces on which mappings
and functionals are defined is dropped and the consideration is
purely algebraic below. Homogeneous idempotent functionals are also
considered. They may be with values in semirings or quasirings
noncommutative or non-associative relative to the addition or the
multiplication.  In Section 2 ordered semirings are described,
propositions and theorems about their construction are proved. Skew
idempotent functionals on them are presented in Section 3. Their
categories are studied. Semirings of functionals and their sequences
are investigated as well.
\par Skew idempotent functionals on semirings may be used for studies of structures
of semirings, their homomorphisms and representations. The main
results are Propositions 2.3, 3.17, 3.19, 3.21, 3.24, 3.27, 3.28,
3.30, 3.34, 3.38, Theorems 2.5, 2.10, 2.11, 3.33, 3.39, 3.40 and
Corollary 2.15.
\par All main results of this paper are obtained for the first time.

\section{Ordered non-associative semirings}
To avoid misunderstandings we first present our definitions.
\par {\bf 1. Definitions.} Let $K$ be a non-void set.
If $K$ is supplied with a binary operation corresponding to a
mapping $\mu : K^2\to K$, then $K$ is called a groupoid.
\par If a binary operation $\mu $ is associative $(ab)c=a(bc)$ for every
$a, b, c\in K$, then $K$ is called a semigroup, where $\mu (a,b)$ is
denoted shortly by $ab$.
\par An element $e=e_{\mu }$ in a groupoid $K$ is called neutral (or unit),
if $eb=be=b$ for each $b\in K$. A semigroup with a neutral element
is called a monoid.
\par An element $b\in K$ in a groupoid with a unit is called left or right
invertible, if there exists a left $b_l^{-1}\in K$ or right inverse
$b_r^{-1}\in K$ respectively, i.e. $b^{-1}_lb=e$ or $bb^{-1}_r=e$
correspondingly. If an element is both left and right invertible,
then it is called invertible.
\par A semigroup with a neutral element in which each element is
invertible is called a group.
\par A groupoid with a neutral element in which each equation $ax=b$ or $xa=b$
has a solution is called a quasigroup.
\par Let $K$ be a set and let two operations $+: K^2\to K$ the addition and
$\times : K^2\to K$ the multiplication be given so that $K$ and
$K\setminus \{ 0 \} $ are monoids (or quasigroups or groupoids with
with neutral elements) relative to $+$ and $\times $ correspondingly
with neutral elements denoted by $e_+ =: 0$ and $e_{\times } =: 1$
so that $a\times 0=0\times a=0$ for each $a\in K$ and either the
left distributivity $a(b+c)=ab+ac$ for every $a, b, c\in K$ or the
right distributivity $(b+c)a=ba+ca$ for every $a, b, c\in K$ is
accomplished, then $K$ is called a semiring (or a quasiring
respectively) with either the left or right distributivity
correspondingly. If it is simultaneously right and left
distributive, then it is called simply a distributive semiring (or a
distributive quasiring respectively). If a type of the
distributivity is not mentioned it will be supposed that the left or
right distributivity is accomplished in a semiring or a quasiring.
\par A semiring (or a quasiring) $K$, which is a group relative to
the addition and $K\setminus \{ 0 \} $ is a group (or a quasigroup)
relative to the multiplication, is called a ring (or a
non-associative ring respectively, i.e. non-associative relative to
the multiplication).
\par A semiring $K$ (or a quasiring, or a ring, or a non-associative
ring) having also a structure of a linear space over a field $\bf F$
and such that $\alpha (a+b) = \alpha a+\alpha b$, $1a=a$, $\alpha
(ab)= (\alpha a)b = a(\alpha b)$ and $(\alpha \beta )a=\alpha (\beta
a)$ for each $\alpha , \beta \in \bf F$ and $a, b\in K$ is called a
semialgebra (or a quasialgebra, or and algebra or a non-associative
algebra correspondingly).
\par A set $K$ with binary operations $\mu _1,...,\mu _n$ will also be called an algebraic
object. An algebraic object is commutative relative to an operation
$\mu _p$ if $\mu _p(a,b)=\mu _p(b,a)$ for each $a, b \in K$.
\par A set $K$ is called directed by a relation $\le $, if
it satisfies the following conditions:
\par $(D1)$ if $x\le y$ and $y\le z$, then $x\le z$;
\par $(D2)$ for each $x\in K$ one has $x\le x$;
\par $(D3)$ for each $x, y\in K$ an element $z\in K$ exists so that
$x\le z$ and $y\le z$. \par If $x\le y$ and $x\ne y$, then one
traditionally writes $x<y$.
\par If  a relation $<$ on $K$ satisfies the following conditions:
\par $(LO1)$ if $x<y$ and $y<z$, then $x<z$;
\par $(LO2)$ if $x<y$, then $y<x$ takes no place;
\par $(LO3)$ if $x\ne y$, then either $x<y$ or $y<x$,
\par then $K$ is called linearly ordered, $<$ is a linear ordering on $K$.
\par $(WO).$ If a set $K$ is linearly ordered and each non-void subset $A$ has a least element in $K$
(i.e. $K$-least element), then $K$ is called well-ordered.
\par An algebraic object $K$ with binary operations $\mu _1,...,\mu _n$ is called either directed
or linearly ordered or well-ordered if it is such as a set
correspondingly and its binary operations preserve an ordering: $\mu
_p(a,b)\le \mu _p(c,d)$ for each $p=1,...,n$ and for every $a, b, c
,d\in K$ so that $a\le c$ and $b\le d$ when $a, b, c, d$ belong to
the same linearly ordered set $Z$ in $K$.
\par Henceforward, we suppose that the minimal element in an ordered
$K$ is zero.
\par Henceforth, for semialgebras, non-associative algebras or
quasialgebras $A$ speaking about ordering on them we mean that only
their non-negative cones $K= \{ y: ~ y\in A, 0\le y \} $ are
considered. For non-negative cones $K$ in semialgebras,
non-associative algebras or quasialgebras only the case over the
real field will be considered, since it is the unique topologically
complete field in which the addition and multiplication operations
are compatible with the natural linear ordering, i.e. $K$ is over
$[0,\infty ) = \{ x: ~ 0\le x, x\in {\bf R} \} $.

\par {\bf 2. Definition.} Suppose that a family $\{ K_j: ~ j\in J \} $ of algebraic
objects $K_j$ is given with binary operations $\mu _1,...,\mu _n$,
where $J$ is a directed infinite set, so that algebraic embeddings
\par $(1)$ $t^j_k: K_j\hookrightarrow K_k$ exist for each $k<j$ in
$J$ with
\par $(2)$ $\mu _l(t^j_k(a),t^j_k(b)) = t^j_k( \mu _l(a,b))$ for all
$a, b\in K_j$ and each $l=1,...,n$;
\par $(3)$ $\psi _p: J\to J$ is a surjective bijective monotone
increasing mapping so that $\psi _p(j)\le j$ and $\psi _p(k)< \psi
_p (j)$ for each $k< j\in J$ for all $p=1,...,n$;
\par $(4)$ $\phi _p: J\to J$ is a bijective monotone increasing mapping so that
$j\le \phi _p(j)$ and $\phi _p(k)< \phi _p(j)$ for each $k< j\in J$
and for all $p=1,...,n$. We define a new algebraic object
\par $(5)$ $s(K_j: ~ j\in J)$ elements of which are $y=(y_j: ~
y_j\in K_j ~ \forall j\in J )\in \prod_{j\in J} K_j$ supplied with
binary operations
\par $(6)$ $\mu _p(y,z) = q$ such that $q_{\psi _p(j)} = \mu _p(y_j,t^{\phi _p(j)}_j(z_{\phi
_p(j)}))$ for each $j\in J$, for all $y, z\in s(K_j: ~ j\in J)$ and
each $p=1,...,n$.

\par {\bf 3. Proposition.} {\it Let $\{ K_j: ~ j\in J \} $ be a family of algebraic
objects with operations $\mu _1,...,\mu _n$ and let an algebraic
object $s(K_j: ~ j\in J)$ be as in Definition 2. If each $K_j$ is
directed so that algebraic embeddings $t^j_k: K_j\hookrightarrow
K_k$ are monotone increasing
\par
$(1)$ $t^j_k(a)< t^j_k(b)$ for each $a<b\in K_j$ and all $k<j\in J$ \\
and $J$ is directed, then $s(K_j: ~ j\in J)$ is naturally directed.}
 \par {\bf Proof.} The algebraic object $s(K_j: ~ j\in J)$ can be partially ordered:
 \par $(2)$ $y\le z \in s(K_j: ~ j\in J)$ if and only if $y_j\le z_j\in K_j$ for each $j\in J$.
Each algebraic object $K_j$ is directed, hence an element $u\in
s(K_j: ~ j\in J)$ with $y_j\le u_j$ and $z_j\le u_j$ for each $j\in
J$ satisfies according to Condition $(2)$ the inequalities $y\le u$
and $z\le u$, consequently, the set $s(K_j: ~ j\in J)$ is directed.
\par If $a\le c$ and $b\le d$ in $s(K_j: ~ j\in J)$, then \par $(3)$ $\mu
_p(a_k,t^{\phi _p(k)}_k(b_{\phi _p(k)})) \le \mu _p(c_k,t^{\phi
_p(k)}_k(d_{\phi _p(k)}))$ for each $k\in J$ and each $p=1,...,n$,
since $a_k\le c_k$ and $b_{\phi _p(k)}\le d_{\phi _p(k)}$. Thus $\mu
_p(a,b)\le \mu _p(c,d)$ for each $p=1,...,n$ and hence the algebraic
object $s(K_j: ~ j\in J)$ is directed.

\par {\bf 4. Proposition.} {\it If $\{ K_j: ~ j\in J \} $ is a family of algebraic
objects with operations $\mu _1,...,\mu _n$ and an algebraic object
$s(K_j: ~ j\in J)$ is as in Definition 2 and if either each $K_j$ is
\par $(1)$ right $\mu _k(\mu _l(b,c),a)=\mu _l(\mu _k(b,a),\mu _k(c,a))$ distributive or
\par $(2)$ left $\mu _k(a,\mu _l(b,c)) = \mu _l(\mu _k(a,b),\mu _k(a,c))$
distributive for each $a, b, c\in K_j$ for a pair of binary
operations $(\mu _k, \mu _l)$ with $k\ne l$, where $\psi _l=id$ and
$\phi _l=id$, then $s(K_j: ~ j\in J)$ is right or left distributive
for $(\mu _k, \mu _l)$ respectively.}
\par {\bf Proof.} Since the mappings $\psi _l=id$ and
$\phi _l=id$ are the identities on $J$, $id(j)=j$ for each $j\in J$,
then for arbitrary elements $a, b, c \in s(K_j: ~ j\in J)$ in the
case $(1)$ we get \par $\mu _k(\mu _l(b_j,c_j),t^{\phi
_k(j)}_j(a_{\phi _k(j)}))=\mu _l(\mu _k(b_j,t^{\phi _k(j)}_j(a_{\phi
_k(j)})),\mu _k(c_j,t^{\phi _k(j)}_j(a_{\phi _k(j)})))$, \\
consequently, $\mu _k(\mu _l(b,c),a)=\mu _l(\mu _k(b,a),\mu
_k(c,a))$. Analogously in the case $(2)$ the equalities are
satisfied: \par $\mu _k(a_j, \mu _l(t^{\phi _k(j)}_j(b_{\phi
_k(j)}),t^{\phi _k(j)}_j(c_{\phi _k(j)}))) = \mu _l(\mu
_k(a_j,t^{\phi _k(j)}_j(b_{\phi _k(j)})),\mu _k(a_j,t^{\phi
_k(j)}_j(c_{\phi _k(j)})))$ \\ and hence $\mu _k(a,\mu _l(b,c)) =
\mu _l(\mu _k(a,b),\mu _k(a,c))$.

\par {\bf 5. Theorem.} {\it If each $K_j$ is non-trivial (contains not only
neutral elements relative to binary
operations $\mu _1,...,\mu _n$) and a number $p$ exists, $1\le p\le
n$, such that $j<\phi _p(j)$ for each $j\in J$, then an algebraic
object $s(K_j: ~ j\in J)$ from Definition 2 is non-associative
relative to a binary operation $\mu _p$.}
\par {\bf Proof.} For arbitrary three elements $a, b, c \in s(K_j: ~ j\in J)$
we deduce from Definition 2:
\par $(1)$ $\mu _p(\mu _p(a,b),c)) =q$ with
\par $(2)$ $q_{\psi _p(\psi _p(j))} = \mu _p(\mu _p(a_j,t^{\phi _p(j)}_j(b_{\phi
_p(j)})),t^{\phi _p(\psi _p(j))}_{\psi _p(j)}(c_{\phi _p(\psi _p
(j))})))\in K_{\psi _p(\psi _p(j))}$ and
\par $(3)$ $\mu _p(a,\mu _p(b,c))=v$ with \par $(4)$ $v_{\psi _p(\psi _p(j))} =
\mu _p (a_{\psi _p(j)},\mu _p(t^{\phi _p(\psi _p(j))}_{\psi _p
(j)}(b_{\phi _p(\psi _p(j))}),t^{\phi _p(\phi _p(\psi _p(j)))}_{\phi
_p (\psi _p(j))}(c_{\phi _p(\phi _p(\psi _p(j)))})))\in K_{\psi _p
(\psi _p(j))}$.
\par The associativity of $\mu _p$, i.e. the equality
$q=v$ is equivalent to \par $(5)$ $[\forall j\in J ~ ~
q_{\psi _p(\psi _p(j))} = v_{\psi _p(\psi _p(j))}]$. \\
Elements $a, b, c$ are arbitrary and generally non-equal and
independent. Therefore, Formulas $(1-5)$ imply that $a_j=a_{\psi
_p(j)}$ for each $j\in J$, consequently, $\psi _p(j)=j$ for each
$j\in J$. Then we get $c_{\phi _p(j)}=c_{\phi _p (\phi _p (j))}$ for
each $j\in J$, but $k< \phi _p(k)$ for each $k\in J$ by the
conditions of this theorem and generally there exists $c_l\ne c_k$
for $l\ne k\in J$, in particular, with $l=\phi _p(j)$ and $k=\phi
_p(\phi _p(j))$. Thus Condition $(5)$ is not satisfied,
consequently, the algebraic object $s(K_j: ~ j\in J)$ is
non-associative relative to the binary operation $\mu _p$.

\par {\bf 6. Theorem.} {\it If $ \{ X_k: ~ k \in K \} $ is a family
of non-associative non-trivial algebraic objects $X_k$ with binary
operations $\mu _1,...,\mu _n$, $~K$ is a set, then there exists a
non-associative algebraic abject $K$ with binary operations $\mu
_1,...,\mu _n$ so that $K$ is not isomorphic with $X_k$ for each
$k\in K$ even up to a direct sum decomposition with either some
associative algebraic objects or the product $X_k^m$, $m\in \bf N$,
with binary operations $\mu _1,...,\mu _n$.}
\par {\bf Proof.} For the set $P := \bigcup_{k\in K} X_k$ one can take
a set $M$ such that $|P|< |M|$ and $\aleph _0 < |M|$, where $|M|$
denotes the cardinality of $M$. Therefore, the inequality
$|P|<|H\times M|$ is satisfied, where $H=\bigcup_{l\in {\bf N}}
K^l$, $~H\times M$ denotes the cartesian product of the sets $H$ and
$M$, $~K^l$ denotes the $l$-fold cartesian product of $K$ with
itself.
\par Let $H$ be directed by inclusion $u\le v$ if and only if $u\subset v$ for any $u, v \in H$.
Each set can be well-ordered in accordance with Zermelo's theorem.
Therefore, let $M$ be well-ordered. The cartesian product $J :=
H\times M$ can be directed lexicographically: $(u,l) < (v,k)$ if and
only if either $u<v$ or $u=v$ and $l<k$, where $l, k \in M$ and $u,
v\in H$. Without loss of generality we can suppose that $J$ has not
a maximal element. Then we put $K_j=\bigoplus_{i\in u} X_i$ for each
$j=(u,l) \in J$ and take the algebraic object $s(K_j: ~ j\in J)$.
Each algebraic object $X_k$ is non-trivial by the conditions of this
theorem. We take $\phi _p$ on $J$ such that $j<\phi _p(j)$ for each
$j\in J$ for all $p=1,...,n$. The cardinality of $s(K_j: ~ j\in J)$
is greater than the cardinality of each $X_k$, consequently, $s(K_j:
~ j\in J)$ is not isomorphic to $X_k$ for each $k\in K$.\par On the
other hand, $s(K_j: ~ j\in J)$ is non-associative relative to each
operation $\mu _1,...,\mu _n$ by Theorem 5. In view of Formulas
5$(1-5)$ even if $s(K_j: ~ j\in J)$ contains either some associative
sub-object $A$ or the product $A=X_k^m$ having a direct sum
complement $B$ so that $s(K_j: ~ j\in J)=A\oplus B$, then $|A|\le
|B|$ and hence $|B|=|s(K_j: ~ j\in J)|$, since $\aleph _0<|M|\le
|H\times M|$. Therefore, the inequality $|X_k|< |B|$ is satisfied
for each $k\in K$, hence $B$ is not isomorphic with $X_k$ for each
$k\in K$.

\par {\bf 7. Proposition.} {\it If $ \{ K_j: ~ j\in J \} $ is a family
of linearly or well-ordered algebraic objects with binary operations
$\mu _1,...,\mu _n$ preserving ordering, i.e. $\mu _p(a,b)<\mu
_p(c,d)$ for each $p=1,...,n$ when either $a<c$ and $b\le d$ or
$a\le c$ and $b<d$, where $J\subset \bf Z$ is a countable set, then
their product $K = \prod_{j\in J} K_j$ is the algebraic object with
binary operations $\mu _1,..,\mu _n$ and it can be naturally
well-ordered.}
\par {\bf Proof.} The product $K$ can be supplied with the binary
operations \par $(1)$ $\mu _p(a,b) = (\mu _p(a_j,b_j)\in K_j: ~ j
\in J)$ for each $p=1,...,n$, for every $a, b\in K$, $ ~ a=(a_j: ~
\forall j ~ a_j\in K_j)$.
\par By the Zermelo's theorem a set $J$ can be well-ordered. Since
$J$ is well-ordered, there exists a $J$-least element $j=\inf
B_{y,z}$ of the subset $B_{y,z} := \{ l: ~ l\in J, ~ y_l\ne z_l \}
$, where $y, z\in s(K_j: ~ j\in J)$ are two arbitrary chosen
elements, $y=(y_l: ~ \forall l\in J ~ y_l\in K_l)$. On the other
hand, $J\subset \bf Z$ by the conditions of this proposition,
consequently, $j\in B_{y,z}$. There exists the lexicographic linear
ordering:
 \par $(2)$ $y<z \in s(K_j: ~ j\in J)$ if there exists $j(y,z)=j\in J$ such
 that $y_k=z_k$ for each $k<j$, whilst $y_j<z_j$.
 As usually we put $y\le z$ if $y=z$ or $y<z$. The set $K$ is
 well-ordered, if $K_j$ is well-ordered for each $j\in J$, since if $A\subset K$ and
 $B=\{ k\in J: ~\exists y, z\in A; ~ y_k\ne z_k \} $, then $j=\inf B\in J$ and
 $\inf A_j\in K_j$, where $A_j = \{ y_j\in K_j: ~y\in A \} $.
\par Let either $a<c$ and $b=d$ or $a=c$ and $b<d$, then from Formulas $(1,2)$ and
the conditions of this proposition it follows that $\mu _p(a,b)< \mu
_p(c,d)$ in $K$. If $a<c$ and $b<d$ and $j(a,c)\le j(b,d)$, then
$a_l<c_l$ and $b_l\le d_l$ for $l=j(a,c)$, consequently, $\mu
_p(a_l,b_l)<\mu _p(c_l,d_l)$. Analogously, if $a< c$ and $b<d$ and
$j(b,d)\le j(a,c)$, then $a_l\le c_l$ and $b_l< d_l$ for $l=
j(b,d)$, consequently, $\mu _p(a_l,b_l)<\mu _p(c_l,d_l)$. Thus the
binary operation $\mu _p$ preserves the ordering on $K$ for each
$p=1,...,n$.

\par {\bf 8. Definition.} A subalgebra $A$ (or a subquasiring, or a subsemiring) of an (a non-associative)
algebra $K$ (or a quasiring, or a semiring) is called a left or
right ideal if $AK\subseteq A$ or $KA\subseteq A$ respectively,
where $AB := \{ c\in K: ~ c=ab, a\in A, b\in B \} $. If $A$ is a
left and right ideal simultaneously, then $A$ is called an ideal.
\par An (a non-associative) algebra (or a quasiring, or a semiring) is called simple if it does not contain ideals
different from $\{ 0 \} $ and $K$, where $K$ is non-trivial, $K\ne
\{ 0 \} $.
\par {\bf 9. Definitions.} Let $s_a(K_j:
~ j\in J )$ denote the subalgebraic object consisting of all
elements $y\in s(K_j: ~ j\in J )$ (see Definition 2) such that
$y=(y_j: ~ \forall j\in J ~ y_j\in K_j)$ with a finite set $J_y :=
\{ j: ~ j\in J, y_j\ne 0 \} $, where $0 = e_{\mu _1,j}\in K_j$ is a
neutral element in $K_j$ relative to a binary operation $\mu _1$ for
each $j$. \par One says that a quasiring or an algebraic object $K$
or a non-associative algebra over a field $\bf F$ is not finitely
(countably) generated, if for each finite (countable respectively)
family of its elements $a_1, a_2,... \in K$ a minimal subquasiring
or a minimal subalgebraic object $Z$ in $K$ or a subalgebra $Z$ over
the same field $\bf F$ correspondingly containing these elements,
$a_1, a_2,... \in Z$, does not coincide with $K$, $Z\ne K$.

\par {\bf 10. Theorem.} {\it If $X$ is a simple non-associative
algebra or a quasiring so that $X$ is not finitely generated, then
there exist a simple non-associative algebra or a quasiring $K$
correspondingly and an embedding $X\hookrightarrow K$ so that $X$ is
not isomorphic with $K$.}
\par {\bf Proof.} We take any infinite set $I$ such that $\aleph _0\le |I|<
|2^X| \aleph _0$. Without loss of generality we suppose that a set
$I$ is well-ordered, since by Zorn's theorem it can be well-ordered.
Let $X$ contain a neutral element $0$ by the addition $+=\mu _1$.
Otherwise one can enlarge $X$ supplying it with a neutral element.
Therefore, without loss of generality $0\in X$. Then we choose a
subset $J\subset I$ without a maximal element (and without a minimal
element if necessary) so that $|J|=|I|$. Then we consider $K :=
s_a(X_j: ~ j\in J ) $, where $X_j=X$ for each $j\in J$ (see
Definition 9). For an addition operation $\mu _1$ we take $\psi
_1(j)=j$ and $\phi _1(j)=j$ for each $j\in J$, while for a
multiplication operation $\mu _2$ we choose $\psi _2(j)\le j$ and
$j<\phi _2(j)$ for each $j\in J$ (see Definition 2). Since $X_j=X$
for each $j\in J$ we put $t^j_k=id$ for every $k\le j\in J$. \par
For each natural number $n$ any family of elements $a_1,...,a_n$
does not generate $X$ as either the non-associative algebra or the
quasiring or the semiring respectively, since $X$ is not finitely
generated.
\par In accordance with Theorem 5 the algebraic object $K$ is
non-associative relative to the multiplication $\mu _2$. By
Proposition 4 $K$ is a non-associative algebra or a quasiring
correspondingly.
\par From Definition 2 it follows that
$\theta (\mu _p(x,z)) = \mu _p(\theta (x),\theta (z))$ for each $x,
z\in X$ for $p=1$ and $p=2$, where \par $(1)$ $\theta (x) = (x_j: ~
x_j=x \forall j\in J)$, $\theta : X\hookrightarrow K$. On the other
hand, if $X$ is a non-associative algebra, then $\theta (\alpha x) =
(y_j: ~ \forall j\in J ~~ y_j=\alpha x ) = \alpha \theta (x)$  for
each $\alpha \in \bf F$ and $x\in X$. Thus $\theta $ is the
embedding of the non-associative algebra or a quasiring $X$ into the
non-associative algebra or a quasiring $K$ correspondingly.
\par Let $M$ be a non-trivial ideal in $K$, that is
$M\ne \{ 0 \} $. We take a family \par $(2)$ $T$ of all elements
$y\in K$ of the form $y=(y_j: ~ y_j= x\delta_{k,j}, ~ j\in J, k\in
P)$, where $x\in X$, $~\delta _{k,j}=0$ for each $k\ne j$, while
$\delta _{j,j}=1$, $~P$ is a finite subset in $J$, particularly, $P$
may be a singleton. An ideal $M$ is a subalgebra or a subquasiring
correspondingly in $K$. From the inclusions $MK\subseteq M$ and
$KM\subseteq M$ it follows that the products $\mu _2(y,a)$ and $\mu
_2(a,y)$ are in $M$ for each $a\in M$ and $y\in T$. Therefore from
formulas for $\mu _1$ and $\mu _2$ on $K$ it follows that $M$
contains $T$. Thus, the inclusion
\par $(3)$ $\theta (X)\subset M$ is fulfilled.
\par The non-associative algebra or the quasiring $\theta (X)$ is simple and
from Conditions $(1-3)$ it follows by multiplications $\mu _2$ and
additions $\mu _1$, that for each finite subset $P$ in $J$ and every
$x_j\in X$ for each $j\in P$ there exists an element $z\in M$ so
that $z_j=x_j$ for each $j\in P$, consequently, $M=K$. Thus the
non-associative algebra or the quasiring $K$ correspondingly is
simple.

\par {\bf 11. Theorem.} {\it Let $\{ X_i: ~ i\in H \} $ be a family of simple
pairwise non-isomorphic either non-associative algebras or
quasirings and let each $X_i$ be not finitely generated, where $H$
is an infinite linearly ordered set and without a maximal element.
Let also for each $i<k$ an embedding $u^k_i: X_k\hookrightarrow X_i$
exist. Then there is a simple non-associative algebra or a quasiring
$K$ correspondingly so that $X_i$ is not isomorphic with $K$ for
each $i\in H$.}
\par {\bf Proof.} For each $X_i$ consider a set $J_i$ without a maximal element and
a minimal element as $J$ for $X$ in Theorem 9. Then we take the set
$J=\bigcup_{i\in H} (i,J_i)$ and order it lexicographically:
$(i,j)<(k,l)$ if either $i<k\in H$ or $i=k$ and $j<l\in J_i$. Let
$K_{i,j} = X_i$ for each $j\in J_i$ and let $K=s_a(K_{i,j}: ~
(i,j)\in J)$.
\par We then choose an addition operation $\mu _1$ with $\psi _1(i,j)=(i,j)$ and $\phi
_1(i,j)=(i,j)$ for each $(i,j)\in J$, while for a multiplication
operation $\mu _2$ let  $\psi _2(i,j)=(l,m)$ with $l\le i$ and $m\le
j$, whilst $\phi _2(i,j)=(v,w)$ with $i<v$ and $j<w$ for each
$(i,j)\in J$. We put $t^{(k,w)}_{(i,j)}=u^k_i$ for every $(i,j)\le
(k,w)\in J$, where $u^i_i=id$. The rest of the proof is analogous to
that of Theorem 9. It shows that $K$ is simple and $K$ is not
isomorphic with $X_i$ as a non-associative algebra or a quasiring
respectively for each $i\in H$.
\par {\bf 12. Corollary.} {\it A family of all simple
either non-associative algebras ${\cal F}_a$ or quasirings ${\cal
F}_q$, which are not pairwise isomorphic, is not a set, i.e. ${\cal
F}_a$ and ${\cal F}_q$ are proper classes (in the NBG axiomatic).}
\par {\bf Proof.} This follows from Theorems 5, 6, 10 and
11. Indeed, if $K$ is either a non-associative algebra or a
quasiring, then either a quotient (non-associative) algebra or a
quotient (quasi) ring $K/A$ by its ideal $A$ is an algebra or a
(quasi) ring respectively. If $A$ is a maximal ideal, then $K/A$ is
simple. Even if $K/A$ is associative, there exists a non-associative
algebra or a quasiring $X$ generated by $K/A$ according to Theorem
5. In view of Theorem 10 it is simple and either a non-associative
algebra or a quasiring correspondingly. In view of Theorems 5 and 10
applied by induction there exists a sequence of simple
non-associative algebras or quasirings respectively satisfying
conditions of Theorem 11, since $I$ and hence $J$ can be taken
uncountable for non countably generated $X$. On the other hand, for
a family of simple either non-associative algebras or quasirings
correspondingly linearly ordered by algebraic embeddings a
non-associative algebra or a quasiring respectively non-isomorphic
to any of them exists. That is, the family of all simple either
non-associative algebras or qiasirings is infinite.
\par If $T$ is a set, it can be well-ordered by Zorn's theorem.
Taking a set $J$ in the class $On$ of all ordinals so that $J$ has
the same cardinality as an infinite set $T$ and without a maximal
element, one gets that an infinite set $P$ can be linearly ordered
as $J$ without a maximal element. Therefore, if a family $ \{ K_j: ~
j \in J \} $ of pairwise non-isomorphic non-associative algebras or
quasirings is taken with an infinite set $J$, then either a
non-associative algebra or a quasiring $X$ correspondingly
non-isomorphic to each $K_j$ would exist. Proposition 4.7$(1)$
\cite{mendelsb} states that $|- Ord(X) \supset (X\notin X \& \forall
u (u\in X \supset u\notin u ))$, while Proposition 4.7$(9)$
\cite{mendelsb} asserts: $|- Ord (X)\supset X=On \vee X\in On$. This
gives the contradiction $J\in J$, since $J\in On$. Therefore, $P$
does not belong to any class and hence is not the set (see the
definition in \S 4.1 \cite{mendelsb}). Thus $P$ is not the set in
the NBG axiomatic. That is ${\cal F}_a$ and ${\cal F}_q$ are proper
classes for a family of all pairwise non-isomorphic non-associative
algebras or quasirings respectively.

\par {\bf 13. Remark.} Corollary 12 is consistent with the fact that the class $On$
of all ordinals is not a set in the NBG theory (see Proposition
4.7$(8)$ in \cite{mendelsb}). Thus constructions given above show
that there are many different directed or linearly ordered or well
ordered non-associative algebras and quasirings.
\par The class $On$ of all ordinals has the addition $\mu _1=+_o$ and
the multiplication $\mu _2=\times _o$ operations which are generally
non-commutative, associative, with unit elements $0$ and $1$
respectively, on $On$ the right distributivity is satisfied (see
Propositions 4.29-4.31 and Examples 1-3 in \cite{mendelsb}). \par
For each non-void set $A$ in $On$ there exists $\sup A\in On$ (see
\cite{kunenb}).
\par If $K$ is a linearly ordered non-commutative relative to the addition
semiring (or a quasiring), then the new operation $(a,b)\mapsto \max
(a,b)=: a\bigoplus b$ defines the commutative addition. Then
$c(a\bigoplus b) = \max (ca,cb) = ca \bigoplus cb$ and $(a\bigoplus
b)c = \max (ac,bc) = ac \bigoplus bc$ for every $a, b, c \in K$,
that is $(T,\bigoplus , \times )$ is left and right distributive.

\par {\bf 14. Definitions.} A partial order is a pair $<T,\le >$
such that $T\ne \emptyset $ and $\le $ is a transitive and reflexive
relation on $T$. A pair $<T,\le >$ is a partial order in the strict
sense if and only if it in addition satisfies $\forall p, q (p\le q
\wedge q\le p \to p=q)$. Then one defines $p<q$ if and only if $p\le
q$ and $p\ne q$.
\par A tree is a partial order in the strict sense, such that for
each $x\in T$ the set $\{ y\in T: y<x \} $ is well-ordered.
\par We consider a directed set $K$ which satisfies the condition:
\par $(DW)$ for each linearly ordered subset $A$ in $K$ there exists
a well-ordered subset $B$ in $K$ such that $A\subset B$.

\par {\bf 15. Corollary.} {\it Let $ \{ K_j: ~ j\in J \} $ be a family
of directed algebraic objects satisfying the condition $(DW)$ with
binary operations $\mu _1,...,\mu _n$ so that algebraic embeddings
$t^j_k: K_j\hookrightarrow K_k$ are monotone increasing
\par
$(1)$ $t^j_k(a)< t^j_k(b)$ for each $a<b\in K_j$ and all $k<j\in J$, \\
where $J$ is a directed set. Then the algebraic object $s(K_j: ~
j\in J)$ (see Definition 2) can be directed to satisfy the condition
$(DW)$ as well.}
\par {\bf Proof.} We consider the algebraic object $s(K_j: ~
j\in J)$ directed as in Proposition 3:  $x\le y \in s(K_j: ~ j\in
J)$ if and only if $x_j\le y_j$ for each $j\in J$. Therefore, if
$x\le y$ and $x\ne y$, then there exists $j\in J$ such that $x_j\ne
y_j$, consequently, $x_j<y_j$. We put \par $(2)$ $x<y$ if and only
if $x\le y$ and $x\ne y$.
\par There are natural projections $\pi _j: s(K_j: ~ j\in
J)\to K_j$ so that $\pi _j(y)=y_j$ for each $j\in J$. For a linearly
ordered subset $A$ in $s(K_j: ~ j\in J)\to K_j$ for each projection
$\pi _j(A)=A_j$ one can take a family of all well-ordered subsets
$B_{j,v}\in K_j$ so that $A_j\subset B_{j,v}$ for each $v\in V_j$,
where $V_j$ is a set. Then each set $B_j := \cap_{v\in V_j} B_{j,v}$
is well-ordered and contains $A_j$. Let $C=\prod_{j\in J} B_j$.
Suppose that $E\subset A$, then $g=(g_j: ~ \forall j\in J ~ g_j=\inf
E_j)$ is in $C$, since $g_j\in B_j$ for each $j\in J$. Therefore,
$g\le \inf E$, i.e. $\inf E=:u$ is in $C$, but then there exists
$y\in A$ so that $u\le y$. Put \par $(3)$ $B =  A\cup
(\bigcup_{E\subset A} \inf E)$. \\  If $U\subset B$, then there
exist sets $H$ and $F_b$ for each $b\in H$ so that \par $U = (B\cap
A)\cup (\bigcup_{b\in H} \inf F_b)$, consequently, \par
$(4)$ $\inf U = \inf ((B\cap A)\cup (\bigcup_{b\in H} F_b))\in B$. \\
On the other hand, if $f, h\in B$, then \par $(5)$ $f=\inf F_f$ and
$h=\inf F_h$ for some suitable subsets $F_f$ and $F_h$ in $A$.
Suppose that these sets $F_h$ are chosen to satisfy Conditions
$(4,5)$. There may be several variants. If $F_f\subset F_h$, then
$h\le f$; if $F_h\subset F_f$, then $f\le h$. If $F_f\setminus F_h$
contains and element $y$ so that $y<c$ for each $c\in F_h$, then
$f\le h$. If $F_h\setminus F_f$ contains an element $y$ so that
$y<c$ for each $c\in F_f$, then $h\le f$. If for each $a\in F_h$
there exists $c\in F_f$ such that $c\le a$ and for each $e\in F_f$
there exists $q\in F_h$ such that $q\le e$, then $f=h$. Together
with $(2)$ this gives that either $f<h$ or $h<f$ or $f=h$, i.e. $B$
is linearly ordered and together with $(4)$ this implies that $B$ is
well-ordered and $(3)$ means that $A\subset B$. Thus the algebraic
object $s(K_j: ~ j\in J)$ is directed so that it satisfies the
condition $(DW)$.

\section{Skew idempotent functionals}
\par {\bf 1. Definitions.} Let $K$ be a well-ordered or directed satisfying
condition 2.14$(DW)$ either semiring or quasiring (or a non-negative
cone in a quasialgebra over the real field $\bf R$) such that
\par $(1)$ $\sup E\in K$ for each  $E\in T$, where $T$ is a family
of subsets of $K$.
\par For a set $X$ and a semiring (or quasiring) $K$
let $C(X,K)$ denote a semiring (or a quasiring respectively) of all
mappings $f: X\to K$ with the point-wise addition
$(f+g)(x)=f(x)+g(x)$ and the point-wise multiplication
$(fg)(x)=f(x)g(x)$ operations for every $f, g \in C(X,K)$ and $x\in
X$. \par If $K$ is a directed semiring (or a directed quasiring) and
$X$ is a linearly ordered set, $C_+(X,K)$ (or $C_-(X,K)$) will
denote the set of all monotone non-decreasing (or non-increasing
correspondingly) maps $f\in C(X,K)$.
\par For the space $C(X,K)$ (or $C_+(X,K)$
or $C_-(X,K)$) we suppose that \par $(2)$ a family $T$ of subsets of
$K$ contains the family $ \{ f(X): ~ f \in C(X,K) \} $ (or $ \{
f(X): ~ f \in C_+(X,K) \} $ or $ \{ f(X): ~ f \in C_-(X,K) \} $
correspondingly) and $K$ satisfies Condition $(1)$.
\par Henceforward, we suppose that the minimal element in $K$ is
zero. \par If a non-associative algebra or a quasialgebra $A$ over
$\bf R$ is considered, we suppose that $K$ is a cone of non-negative
elements (non-negative cone) in it (see also \S 2.1).
\par {\bf 2. Remark.} As an example of a semiring (or a quasiring) $K$ in
Definitions 1 one can take $K=On$ or $K = \{ A: ~ A\in On, |A|\le b
\} $, where $b$ is a cardinal number such that $\aleph _0\le b$ (see
also Remark 2.13).  Evidently, $K=On$ satisfies Condition 1$(1)$,
since $\sup E$ exists for each set $E$ in $On$ (see \cite{kunenb}).
\par  Another examples are provided by Theorems 2.5, 2.10, 2.11,
Propositions 2.3, 2.4, 2.7 and Corollaries 2.12, 2.15.
\par It is possible to modify Definition 1 in the following manner.
For a well-ordered $K$ without Condition 1$(1)$ one can take the
family of all bounded functions $f: X\to K$ and denote this family
of functions by $C(X,K)$ for the uniformity of the notation.
\par For a directed $K$ satisfying Condition 2.14$(DW)$ without
Condition 1$(1)$ it is possible to take the family of all monotone
non-decreasing (or non-increasing) bounded functions $f: X\to K$ for
a linearly ordered set $X$ and denote this family by $C_+(X,K)$ (or
$C_-(X,K)$ correspondingly) also.
\par Naturally, $C(X,K)$ has also the structure of the left and
right module over the semiring (or the quasiring correspondingly)
$K$, i.e. $af$ and $fa$ belong to $C(X,K)$ for each $a\in K$ and
$f\in C(X,K)$. To any element $a\in K$ the constant mapping $g^a\in
C(X,K)$ corresponds such that $g^a(x)=a$ for each $x\in X$. If $K$
is right (or left) distributive, then $q(f+h)=qf+qh$ (or
$(f+h)q=fq+hq$ correspondingly) for every $q, f, h\in C(X,K)$.
\par The semiring (or the quasiring) $C(X,K)$ will be considered
directed: \par $(1)$ $f\le g$ if and only if $f(x)\le g(x)$ for each
$x\in X$. \par Indeed, if $f, h \in C(X,K)$, then $a=\sup (f(X))\in
K$ and $b=\sup (h(X))\in K$ according to Condition 1$(1)$. Then
there exists $c\in K$ so that $a\le c$ and $b\le c$, consequently,
$f\le g^c$ and $h\le g^c$. Thus for each $f, h\in C(X,K)$ a function
$q\in C(X,K)$ exists so that $f\le q$ and $h\le q$. From $a+b\le
c+d$ and $ac\le bd$ for each $a\le c$ and $b\le d$ in $K$ it follows
that $f+q\le g+h$ and $fq\le gh$ for each $f\le g$ and $q\le h$ in
$C(X,K)$. \par If $f\le g$ and $f\ne g$ (i.e. $\exists x\in X ~
f(x)\ne g(x)$), then we put $f<g$.
\par For a mapping $f\in C(X,K)$ its support $supp (f)$ is defined as usually
\par $(2)$ $supp (f) := \{ x: ~ x\in X, f(x)\ne 0 \} $.
\par {\bf 3. Lemma.} {\it If $E$ is a subset in $X$, then
$C(X,K|E) := \{ f: ~ f\in C(X,K), supp (f)\subset E \} $ is an ideal
in $C(X,K)$.}
\par {\bf Proof.} If $f\in C(X,K|E)$ and $g\in C(X,K)$, then
$f(x)g(x)=0$ and $g(x)f(x)=0$ when $f(x)=0$, consequently, $supp
(fg)$ and $supp (gf)$ are contained in $E$. Moreover, if $f, h\in
C(X,K|E)$, then $supp (f+h)$ and $supp (h+f)$ are contained in $E$,
since $f(x)+h(x)=0$ and $h(x)+f(x)=0$ for each $x\in X\setminus E$.
Thus $C(X,K|E)$ is a semiring (or a quasiring respectively) and
$C(X,K|E)C(X,K)\subseteq C(X,K|E)$ and $C(X,K)C(X,K|E)\subseteq
C(X,K|E)$.
\par {\bf 4. Corollary.} {\it If $E$ is a subset in $X$, then
$C(E,K)$ is an ideal in $C(X,K)$.}
\par {\bf Proof.} For a subset $E$ in $X$ one gets $C(E,K)$ isomorphic with
$C(X,K|E)$, since each $f\in C(E,K)$ has the zero extension on
$X\setminus E$.
\par {\bf 5. Lemma.} {\it For a linearly ordered set $X$ and a
directed semiring (quasiring) $K$ there are directed semirings (or
quasirings correspondingly) $C_+ (X,K)$ and $C_-(X,K)$.}
\par {\bf Proof.} The sets $C_+(X,K)$ and $C_-(X,K)$ are directed
according to Condition 2$(1)$ with a partial ordering inherited from
$C(X,K)$. Since $a+b\le c+d$ and $ac\le bd$ for each $a\le c$ and
$b\le d$ in $K$, then $f+q\le g+h$ and $fq\le gh$ for each $f\le g$
and $q\le h$ in $C_+(X,K)$ and $C_-(X,K)$. On the other hand, for
each $f, h \in C(X,K)$ there exists $g^c\in C(X,K)$ so that $f\le
g^c$ and $h\le g^c$ (see \S 2). If $f(x)\le f(y)$ and $h(x)\le h(y)$
for $f, h\in C_+(X,K)$ and each $x\le y$ in $X$, then $f(x)+h(x)\le
f(y)+h(y)$ and $f(x)h(x)\le f(y)h(y)$, consequently, $f+h$ and $fh$
are in $C_+(X,K)$. Analogously, if $f, h \in C_-(X,K)$, then $f+h$
and $fh$ are in $C_-(X,K)$.  But a constant mapping $g^c$ belongs to
$C_+(X,K)$ and $C_-(X,K)$. Thus $C_+(X,K)$ and $C_-(X,K)$ are
directed semirings (or quasirings correspondingly).
\par {\bf 6. Lemma.} {\it Suppose that $K$ satisfies Conditions 2.14$(DW)$
and 1$(1,2)$. Then the functions
\par $(1)$ $f\vee g (x) := \max (f(x), g(x))$ and
\par $(2)$ $f\wedge g (x) := \min (f(x), g(x))$ \\
are in $C(X,K)$ (or in $C_+(X,K)$ or in $C_-(X,K)$) for every pair
of functions $f, g \in C(X,K)$ (or in $C_+(X,K)$ or in $C_-(X,K)$
correspondingly) satisfying the condition: \par $(3)$ for each $x\in
X$ either $f(x)<g(x)$ or $g(x)<f(x)$ or $f(x)=g(x)$.}
\par {\bf Proof.} Let $f, g \in C(X,K)$ satisfy
Condition $(3)$. Then the sets $ \{ x: ~ x\in X, f(x)\le g(x) \} $
and $ \{ x: ~ x\in X, f(x)\le g(x) \} $ are contained in $X$. For
each subset $E$ in $K$ the sets
\par $(f\vee g)^{-1}(E) = [f^{-1}(E)\cap \{ x: ~ x\in X, g(x)\le f(x)
\} ] \cup [g^{-1}(E)\cap \{ x: ~ x\in X, f(x)\le g(x) \} ]$ and
\par $(f\wedge g)^{-1}(E) = [f^{-1}(E)\cap \{ x: ~ x\in X, f(x)\le g(x)
\} ] \cup [g^{-1}(E)\cap \{ x: ~ x\in X, g(x)\le f(x) \} ]$
\\ are contained in $X$, consequently, the mappings $f\vee g$ and $f\wedge
g$ are in $C(X,K)$, since either $\sup f \le \sup g$ or $\sup g
<\sup f$ for $f$ and $g$ satisfying $(3)$ and hence $\max (\sup f,
\sup g) \in K$ and $\min (\sup f, \sup g)\in K$.
\par If $f, g\in C_+(X,K)$ and $x<y\in X$, then $f(x)\le f(y)$ and
$g(x)\le g(y)$. If $f(x)\le g(x)$ and $g(y)\le f(y)$, then $f\vee
g(x)=g(x)\le g(y)\le f(y) = f\vee g(y) $ and $f\wedge g(x)=f(x)\le
g(x)\le g(y)=f\wedge g(y)$. If $f(x)\le g(x)$ and $f(y)\le g(y)$,
then $f\vee g(x)=g(x)\le g(y)=f\vee g(y)$ and $f\wedge g(x)=f(x)\le
f(y)=f\wedge g(y)$. Therefore, $f\vee g(x)\le f\vee g(y)$ and
$f\wedge g(x)\le f\wedge g(y)$. Thus $f\vee g$ and $f\wedge g\in
C_+(X,K)$. Analogously if $f, g\in C_-(X,K)$, then $f\vee g$ and
$f\wedge g\in C_-(X,K)$.

\par {\bf 7. Notation.} Let $\bigodot $ denote the mapping on $[K\times
C(X,K)]\cup [C(X,K)\times K]$ with values in $C(X,K)$ such that
\par $(1)$ $c\bigodot f= g^c+f$  and $f\bigodot c= f+ g^c$ for each
$c\in K$ and $f\in C(X,K)$, where $g^c(x)=c$ for each $x\in X$,
whilst the sum is taken point-wise $(f+g)(x)=f(x)+g(x)$ for every
$f, g \in C(X,K)$ and $x\in X$.
\par {\bf 8. Definition.} We call a mapping $\nu $ on $C(X,K)$
(or $C_+(X,K)$ or $C_-(X,K)$) with values in $K$ an idempotent
($K$-valued) functional if it satisfies for each $f, g, g^c \in
C(X,K)$ (or in $C_+(X,K)$ or $C_-(X,K)$ correspondingly) the
following five conditions
\par $(1)$ $\nu (g^c)=c$;
\par $(2)$ $\nu (c\bigodot f) = c\bigodot \nu (f)$ and
\par $(3)$ $\nu (f\bigodot c) = \nu (f)\bigodot c$;
\par $(4)$ $\nu (f\vee g)=\nu (f)\vee \nu (g)$ when $f, g$ satisfy Condition 6$(3)$ and
\par $(5)$ $\nu (f\wedge g)=\nu (f)\wedge \nu (g)$ if $f, g$ satisfy Condition 6$(3)$,\\
where $a\vee b =\max (a,b)$ and $a\wedge b=\min (a,b)$ for each $a,
b\in K$ when either $a<b$ or $a=b$ or $b<a$.
\par A mapping (functional) $\nu $ on $C(X,K)$ (or $C_+(X,K)$ or $C_-(X,K)$)
with values in $K$ we call order preserving (non-decreasing), if
\par $(6)$ $\nu (f)\le \nu (g)$ for each $f\le g$ \\ in
$C(X,K)$ (or $C_+(X,K)$ or $C_-(X,K)$ respectively), i.e. when
$f(x)\le g(x)$ for each $x\in X$. \par A functional $\nu $ is called
left or right $K$-homogeneous on $C(X,K)$ (or $C_+(X,K)$ or
$C_-(X,K)$) if
\par $(7)$ $\nu (bf)=b \nu (f)$ or
\par $(8)$ $\nu (fb)=\nu (f)b$ \\ for each $f$ in $C(X,K)$ (or $C_+(X,K)$
or $C_-(X,K)$ correspondingly) and $b\in K$. A functional left and
right homogeneous simultaneously is called homogeneous.
\par {\bf 9. Remark.} If a functional satisfies Condition
10$(4)$, then it is order preserving.
\par  The Dirac functional $\delta _x$ defined by
the formula:
\par $(1)$ $\delta _xf=f(x)$ \\ is the idempotent $K$-homogeneous
functional on $C(X,K)$, where $x$ is a marked point in $X$.
\par If functionals $\nu _1,...,\nu _n$ are idempotent
and the multiplication in $K$ is distributive, then for each
constants \par $(2)$ $c_1>0,...,c_n>0$ in $K$ with
\par $(3)$ $c_1+...+c_n=1$ functionals
\par $(4)$ $c_1\nu _1+...+c_n\nu _n$ and \par $(5)$ $\nu _1c_1+...+\nu _nc_n$ \\ are
idempotent. Moreover, if the multiplication in $K$ is associative
and distributive and constants satisfy Conditions $(2,3)$ and
functionals $\nu _1,...,\nu _n$ are $K$-homogeneous, then
functionals of the form $(4,5)$ are also $K$-homogeneous.
\par The considered here theory is different from the usual real field $\bf R$,
since $\bf R$ has neither an infimum nor a supremum, i.e. it is not
well-ordered and satisfy neither 2.14$(DW)$ nor 1$(1)$.
\par {\bf 10. Lemma.} {\it Suppose that either \par $(1)$ $K$ is well-ordered
and satisfies Conditions 1$(1,2)$ or \par $(2)$ $X$ is linearly
ordered and $K$ is directed and satisfies Conditions 2.14$(DW)$ and
1$(1,2)$. Then there exists an idempotent $K$-homogeneous functional
$\nu $ on $C(X,K)$ in case $(1)$, on $C_+(X,K)$ and $C_-(X,K)$ in
case $(2)$. Moreover, if $K\in On$ and $K$ is infinite, $\aleph
_0\le |K|$, or $K=On$, $X$ is not a singleton, $|X|>1$, then $\nu $
has not the form either 9$(4)$ or 9$(5)$ with Dirac functionals $\nu
_1,...,\nu _n$ relative to the standard addition in $On$.}
\par {\bf Proof.} Suppose that $\nu $ is an order preserving functional on $C(X,K)$
(or $C_+(X,K)$ or $C_-(X,K)$). If functions $f, g$ in $C(X,K)$ (or
$C_+(X,K)$ or $C_-(X,K)$ respectively) satisfy Condition 6$(3)$,
then in accordance with Lemma 6 there exists $f\vee g$ and $f\wedge
g$ in the corresponding $C(X,K)$ (or $C_+(X,K)$ or $C_-(X,K)$).
Since $f\vee g\ge f$ and $f\vee g\ge g$ and $f\wedge g\le f$ and
$f\wedge g\le g$ and the functional $\nu $ is order preserving, then
$\nu (f)\vee \nu (g) \le \nu (f\vee g)$ and $\nu (f\wedge g)\le \nu
(f)\wedge \nu (g)$.
\par  Let also $E$ be a subset in $X$, we put
\par $(3)$ $\nu (f)=\nu _{E} (f)=\sup_{x\in E} f(x)$.
\\ This functional exists due Conditions 1$(1,2)$, since in both cases
$(1)$ and $(2)$ of this lemma, the image $f(E)$ is linearly ordered
and is contained in $K$. \par From the fact that the addition
preserves ordering on $K$ (see \S 2.1) it follows that Properties
$(1-3,7,8)$ are satisfied for the functional $\nu $ given by Formula
$(3)$. If $f\le g$ on $X$, then for each $a\in f(E)$ there exists
$b\in g(E)$ so that $a\le b$, consequently, $\nu (f)\le \nu (g)$,
i.e 8$(6)$ is fulfilled. \par We consider any pair of functions $f,
g$ in $C(X,K)$ (or $C_+(X,K)$ or $C_-(X,K)$) satisfying Condition
6$(3)$. In case $(2)$ a set $X$ is linearly ordered, in case $(1)$
$K$ is well-ordered, hence $f(X)$, $g(X)$, $f(E)$ and $g(E)$ are
linearly ordered in $K$. Then for each $a\in f(E)\cup g(E)$ there
exist $b \in (f\vee g)(E)$ so that $a\le b$, while for each $c\in
(f\vee g)(E)$ there exists $d\in f(E)\cup g(E)$ so that $c\le d$,
hence $\nu (f\vee g)= \nu (f)\vee \nu (g)$. Moreover, for each $a\in
f(E)\cup g(E)$ there exists $b\in (f\wedge g)(E)$ so that $b\le a$
and for each $c\in (f\wedge g)(E)$ there exists $d\in f(E)\cup g(E)$
so that $d\le c$, consequently, $\nu (f\wedge g)=\nu (f)\wedge \nu
(g)$. Thus Properties 8$(4,5)$ are satisfied as well.
\par If a set $X$ is not a singleton, $|X|>1$, and $K\subset On$ is infinite,
$\aleph _0\le |K|$, then taking a set $E$ in $X$ different from a
singleton, $|E|>1$, we get that the functional given by Formula
$(3)$ can not be presented with the help of Dirac functionals $\nu
_1=\delta_{x_1},...,\nu _n=\delta_{x_n}$ by Formula either 9$(4)$ or
9$(5)$ relative to the standard addition in $On$, since functions
$f$ in $C(X,K)$ (or $C_+(X,K)$ or $C_-(X,K)$) separate points in $X$
(see Remark 2).
\par {\bf 11. Remark.} Relative to the idempotent addition
$x\vee y=\max (x,y)$ the functional $\nu _E$ given by 10$(3)$ has
the form $\nu _E(f) = \vee_{x\in E} \delta_x (f)$. \par Let $I(X,K)$
denote the set of all idempotent $K$-valued functionals on $C(X,K)$,
$~I_l(X,K)$ of all idempotent $K$-valued functionals on $C_+(X,K)$,
let also $I_h(X,K)$ and $I_{l,h}(X,K)$ denote their subsets of
idempotent homogeneous functionals.
\par {\bf 12. Definitions.} A functional $\nu : C(X,K)\to K$ is called weakly additive, if
\par $(1)$ $\nu (h+g^c)=\nu (h)+c$ and  $\nu (g^c+h)=c+\nu (h)$ for all $c\in K$ and $h\in
C(X,K)$;
\par $(2)$ order preserving if $\nu (f)\le \nu (h)$ for each $f\le
h\in C(X,K)$;
\par $(3)$ normalized at $c\in K$, if $\nu (g^c)=c$;
\par $(4)$ non-expanding if $\nu (f)\le \nu (h)+c$ when $f\le
h+g^c$ and $\nu (f)\le c+\nu (h)$ when $f\le g^c+h$ for any $f, h
\in C(X,K)$ and $c\in K$,
\par where $\nu $ may be non-linear as well.
\par The family of all order preserving weakly additive functionals
on $C(X,K)$ (or $C_+(X,K)$) with values in $K$ will be denoted
${\cal O}(X,K)$ (or ${\cal O}_l(X,K)$ respectively).
\par If $E\subset C(X,K)$ (or $E\subset C_+(X,K)$) satisfies the conditions: $g^0\in E$,
$g+b$ and $b+g\in E$ for each $g\in E$ and $b\in K$, then $E$ is
called an ${\cal A}$-subset.
\par {\bf 13. Lemma.} {\it If $\nu : C(X,K)\to K$ is an order preserving weakly
additive functional, then it is non-expanding.}
\par {\bf Proof.} Suppose that $f, h\in C(X,K)$ and $b\in K$ are such
that $f(x)\le (h(x)+c)$ or $f(x)\le (c+h(x))$ for each $x\in X$,
then 12$(1,2)$ imply that $\nu (f)\le (\nu (h)+c)$ or $\nu (f)\le
(c+\nu (h))$ respectively. Thus the functional $\nu $ is
non-expanding.
\par {\bf 14. Lemma.} {\it Suppose that $A$ is an ${\cal A}$-subset (or a left or right
submodule over $K$) in $C(X,K)$ (or in $C_+(X,K))$and $\nu : A\to K$
is an order preserving weakly additive functional (or left or right
$K$-homogeneous with left or right distributive quasi-ring $K$
correspondingly). Then there exists an order preserving weakly
additive (or left or right $K$-homogeneous correspondingly)
functional $\mu : C(X,K)\to K$ (or $\mu : C_+(X,K)\to K$
respectively) such that its restriction on $A$ coincides with $\nu
$.}
\par {\bf Proof.} One can consider the set $\cal F$ of all pairs $(B,\mu )$
so that $B$ is an ${\cal A}$-subset (or a left or right submodule
over $K$ respectively), $A\subseteq B\subseteq C(X,K)$, $\mu $ is an
order preserving weakly additive functional on $B$ the restriction
of which on $A$ coincides with $\nu $. The set $\cal F$ is partially
ordered: $(B_1,\mu _1)\le (B_2,\mu _2)$ if $B_1\subseteq B_2$ and
$\mu _2$ is an extension of $\mu _1$. In accordance with Zorn's
lemma a maximal element $(E,\mu )$ in ${\cal F}$ exists. \par If
$E\ne C(X,K)$, there exists $g\in C(X,K)\setminus E$. Let $E_- := \{
f : ~ f\in E, f\le g \} $ and $E_+ := \{ f : ~ f\in E, g\le f \} $,
then $\mu (h) \le \mu (q)$ for each $h\in E_-$ and $q\in E_+$,
consequently, an element $b\in K$ exists such that $\mu (E_-)\le b
\le \mu (E_+)$. Then we put $F = E\cup \{ g+g^c, g^c+g: c\in K \} $
( or $F$ is a minimal left or right module over $K$ containing $E$
and $g$ correspondingly). Then one can put $\mu (g+g^c)=b+c$ and
$\mu (g^c+g)=c+b$ (moreover, $\mu (d(g+g^c)) = d \mu (g) + dc$ or
$\mu ((g+g^c)d) = \mu (g)d + cd$ for each $d\in K$ correspondingly)
for each $c\in K$. Then $\mu $ is an order preserving weakly
additive functional (left or right homogeneous correspondingly) on
$F$. This contradicts the maximality of $A$. \par For $C_+(X,K)$ the
proof is analogous.
\par {\bf 15. Definitions.} If sets $X$ and $Y$ are
given and $f: X\to Y$ is a mapping, $K_1, K_2 $ are ordered
quasirings (or may be particularly semirings) with an
order-preserving algebraic homomorphism $u: K_1\to K_2$  then it
induces the mapping ${\cal O}(f,u): {\cal O}(X,K_1)\to {\cal
O}(Y,K_2)$ according to the formula: \par $(1)$ $({\cal O}(f,u)(\nu
))(g) = u[\nu (g_1(f))]$ for each $g_1\in C(Y,K_1)$ and $\nu \in
{\cal O}(X,K_1)$, where $u\circ g_1=g\in C(Y,K_2)$, $g_1\in
C(Y,K_1)$, $({\cal O}(f,u)(\nu ))$ is defined on $({\hat f},{\hat
u})(C(X,K_1)) = \{ t: ~ t\in C(Y,K_2); \forall x\in X ~
t(x)=u(h\circ f(x)), h\in C(Y,K_1) \} $.
\par By $I(f,u)$ will be denoted the restriction of ${\cal O}(f,u)$ onto
$I(X,K)$. If $K$ is fixed, i.e. $u=id$, then we write for short
${\cal O}(f)$ and $I(f)$ omitting $u=id$. If $X=Y$ and $f=id$ we
write for short ${\cal O}_2(u)$ and $I_2(u)$ respectively omitting
$f=id$.
\par Let ${\cal S}$ denote a category such that a family $Ob({\cal S})$ of its objects
consists of all sets, a family of morphisms $Mor(X,Y)$ consists of
all mappings $f: X\to Y$ for every $X, Y\in Ob({\cal S})$. \par Let
$\cal K$ be the category objects of which $Ob({\cal K})$ are all
ordered quasirings satisfying Conditions 2.14, $Mor (A,B)$ consists
of all order-preserving algebraic homomorphisms for each $A, B\in
{\cal K}$. Then by ${\cal K}_w$ we denote its subcategory of
well-ordered quasirings and their order-preserving algebraic
homomorphisms.
\par We denote by ${\cal OK}$ a category with the families of objects $Ob({\cal OK}) =
\{ {\cal O}(X,K): X\in Ob({\cal S}), K\in Ob ({\cal K}_w) \} $ and
morphisms $Mor ({\cal O}(X,K_1), {\cal O}(Y,K_2) )$ for every $X,
Y\in Ob({\cal S})$ and $K_1, K_2\in Ob({\cal K}_w)$. Then ${\cal
IK}$ denotes a category with families of objects $Ob({\cal IK}) = \{
I(X,K): X\in Ob({\cal S}), K\in Ob ({\cal K}_w) \} $ and morphisms
$Mor (I(X,K_1), I(Y,K_2) )$ for every $X, Y\in Ob({\cal S})$ and
$K_1, K_2\in Ob({\cal K})$.

\par By ${\cal S}_l$ will be denoted a category objects of which are linearly
ordered sets, $Mor(X,Y)$ consists of all monotone nondecreasing
mappings $f: X\to Y$, that is $f(x)\le f(y)$ for each $x\le y\in X$,
where $X, Y\in Ob({\cal S}_l)$. Then analogously ${\cal O}_l(f,u):
{\cal O}_l(X,K_1)\to {\cal O}_l(Y,K_2)$ for each $X, Y \in Ob({\cal
S}_l)$ and $f\in Mor(X,Y)$, $ ~ K_1, K_2 \in Ob ({\cal K})$, $ ~
u\in Mor (K_1,K_2)$ according to the formula: \par $(2)$ $({\cal
O}_l(f,u)(\nu ))(g) = u[\nu (g_1(f))]$ for each $g_1\in C_+(Y,K_1)$
and $u\circ g_1 = g \in C(Y,K_2)$ and $\nu \in {\cal O}_l(X,K_1)$,
where $({\cal O}_l(f,u)(\nu ))$ is defined on $({\hat f},{\hat
u})(C_+(X,K_1)) := \{ t: ~ t\in C_+(Y,K_2); \forall x\in X ~
t(x)=u(h\circ f(x)), h\in C_+(Y,K_1) \} $. Then the category ${\cal
O}_l{\cal K}$ with families of objects $Ob({\cal O}_l{\cal K}) = \{
{\cal O}_l(X,K): X\in Ob ({\cal S}_l), K\in Ob ({\cal K}) \} $ and
morphisms $Mor ({\cal O}_l(X,K_1), {\cal O}_l(Y,K_2) ) $ and the
category ${\cal I}_l{\cal K}$ with $Ob ({\cal I}_l{\cal K}) = \{
I_l(X,K): X\in Ob ({\cal S}_l), K\in Ob ({\cal K}) \} $ and $Mor
(I_l(X,K_1), I_l(Y,K_2))$ are defined.
\par Subcategories of left homogeneous functionals we denote by
${\cal O}_h{\cal K}$, ${\cal O}_{l,h}{\cal K}$, ${\cal I}_h{\cal
K}$, ${\cal I}_{l,h}{\cal K}$ correspondingly. They are taken on
subcategories ${\cal K}_{w,l}$ in ${\cal K}$ or ${\cal K}_l$ in
${\cal K}$ of left distributive quasirings.

\par {\bf 16. Lemma.} {\it There exist covariant functors ${\cal
O}$, ${\cal O}_h$ and ${\cal O}_l$, ${\cal O}_{l,h}$ in the
categories ${\cal S}$ and ${\cal S}_l$ respectively.}
\par {\bf Proof.} If $X, Y\in Ob({\cal S})$ and $f\in Mor(X,Y)$, $g\le h$ in $C(Y,K)$,
where $K\in Ob ({\cal K}_w)$ (or in ${\cal K}_{w,l}$) is marked,
then $g\circ f \le h\circ f$ in $C(X,Y)$, consequently, $({\cal
O}(f)(\nu ))(g)=\nu (g\circ f) \le \nu (h\circ f)= ({\cal O}(f)(\nu
))(h)$ for each $\nu \in {\cal O}(X,K)$. If $c\in K$, $g^c\in
C(Y,K)$, then $g^c\circ f\in C(X,K)$, also $({\cal O}(f)(\nu ))
(g^c+h) =\nu (g^c\circ f + h\circ f) = c+ \nu (h\circ f) = c+({\cal
O}(f)(\nu ))(h)$ and $({\cal O}(f)(\nu )) (h+g^c) =\nu (h\circ f +
g^c\circ f) = \nu (h\circ f) +c = ({\cal O}(f)(\nu ))(h) +c $ for
each $h\in C(Y,K)$. If $1_X\in Mor (X,X)$, $1_X(x)=x$ for each $x\in
X$, then $1_X\circ q=q$ for each $q\in Mor (Y,X)$ and $t\circ 1_X=t$
for each $t\in Mor (X,Y)$. On the other hand, $({\cal O}(1_X)(\nu
))(g)= \nu (g\circ 1_X)=\nu (g)$ for each $g\in C(X,K)$, i.e. ${\cal
O}(1_X) = 1_{{\cal O}(X)}$. Evidently, $({\cal O}(f\circ s) (\nu
))(g) = \nu (g\circ f\circ s) = ({\cal O}(s)(\nu )(g\circ f) =
(({\cal O}(f)\circ {\cal O}(s))(\nu ))(g)$.
\par If $\nu \in {\cal O}_h(X,K)$, then $({\cal O}(f)(\nu )) (bg) =
\nu (bg\circ f) = b\nu (g\circ f) = (b({\cal O}(f)(\nu ))(g)$. For
${\cal O}_l$ (or ${\cal O}_{l,h}$) the proof is analogous with $X, Y
\in Ob ({\cal S}_l)$, $C_+(X,K)$ and $C_+(Y,K)$, where $K\in Ob
({\cal K})$ (or $K\in Ob ({\cal K}_l)$) is marked.

\par {\bf 17. Proposition.} {\it If $f\in Mor(X,Y)$ for $X, Y\in
Ob({\cal S})$ or in $Ob({\cal S}_l)$, then \par ${\cal
O}(f)(I(X,K))\subseteq I(Y,K)$ and ${\cal O}_h(f)(I_h(X,K))\subseteq
I_h(Y,K)$ for $K\in Ob ({\cal K}_{w,l})$ or ${\cal
O}_l(f)(I_h(X,K))\subseteq I_l(Y,K)$ or ${\cal
O}_{l,h}(f)(I_{l,h}(X,K))\subseteq I_{l,h}(Y,K)$ for $K\in Ob ({\cal
K})$ or $K\in Ob ({\cal K}_l)$ correspondingly.}
\par {\bf Proof.} If $g, h \in C(Y,K)$ are such that $g\vee h$ or $g\wedge h$ exists
(see Condition $(3)$ in Lemma 6) and $f: X\to Y$ is a mapping, $\nu
\in I(X,K)$ (or $I_l(X,K)$), then
\par $({\cal O}(f)(\nu )) (g\vee h) = \nu (g\circ f \vee h\circ f) =
\nu (g\circ f )\vee \nu (h\circ f) = ({\cal O}(f)(\nu )) (g)\vee
({\cal O}(f)(\nu )) (h)$ or \par $({\cal O}(f)(\nu )) (g\wedge h) =
\nu (g\circ f \wedge h\circ f) = \nu (g\circ f )\wedge \nu (h\circ
f) = ({\cal O}(f)(\nu )) (g)\wedge ({\cal O}(f)(\nu )) (h)$. Then
for each $c\in K$  one gets \par $({\cal O}(f)(\nu )) (g^c\bigodot
h) = \nu (g^c\circ f \bigodot h\circ f) = \nu (g^c\circ f )\bigodot
\nu (h\circ f) = c \bigodot ({\cal O}(f)(\nu )) (h)$ and
\par $({\cal O}(f)(\nu )) (h\bigodot g^c) = \nu (h\circ f
\bigodot g^c\circ f) = \nu (h\circ f)\bigodot \nu (g^c\circ f ) =
({\cal O}(f)(\nu )) (h) \bigodot c$.
\par If $\nu \in I_h(X,K)$ (or $I_{l,h}(X,K)$), then $({\cal
O}(f)(\nu ))(bg) = \nu (bg\circ f) = b \nu (g\circ f) = (b({\cal
O}(f)(\nu )))(g)$.

\par {\bf 18. Definitions.} A covariant functor $F: {\cal S}\to {\cal S}$ will be called
epimorphic (monomorphic) if it preserves epimorphisms
(monomorphisms). If $\phi : A\hookrightarrow X$ is an embedding,
then $F(A)$ will be identified with $F(\phi )(F(A))$.
\par If for each $f\in Mor(X,Y)$ and each subset $A$ in $Y$,
the equality $(F(f)^{-1})(F(A)) = F(f^{-1}(A))$ is satisfied, then a
covariant functor $F$ is called preimage-preserving. When
$F(\bigcap_{j\in J} X_j)= \bigcap_{j\in J} F(X_j)$ for each family $
\{ X_j: j\in J \} $ of subsets in $X\in Ob({\cal S})$ (or in $Ob
({\cal S}_l)$), the monomorphic functor $F$ is called
intersection-preserving.
\par If a functor $F$ preserves inverse mapping system limits, it is called continuous.
\par A functor is said to be semi-normal when it is
monomorphic, epimorphic, also preserves intersections, preimages and
the empty space.
\par If a functor is monomorphic, epimorphic, also preserves
intersections and the empty space, then it is called weakly
semi-normal.

\par {\bf 19. Proposition.} {\it The functor ${\cal O}$ (or ${\cal O}_h$,
${\cal O}_l$, ${\cal O}_{l,h}$) is monomorphic.}
\par {\bf Proof.} Consider $X, Y\in Ob({\cal S})$ (or in $Ob ({\cal S}_l)$
respectively) with an embedding $s: X\hookrightarrow Y$
(order-preserving respectively). Suppose that $\nu _1\ne \nu _2\in
{\cal O}(X,K)$ (or in ${\cal O}_h(X,K)$, ${\cal O}_l(X,K)$, ${\cal
O}_{l,h}(X,K)$ correspondingly). This means that a mapping $g\in
C(X,K)$ (or in $C_+(X,K)$ correspondingly)  exists such that $\nu
_1(g)\ne \nu _2(g)$. A function $u\in C(Y,K)$ (or in $C_+(Y,K)$
respectively) exists such that $u\circ s = g$, hence $({\cal
O}(s)(\nu _k)))(u) = \nu _k(u\circ s) = \nu _k(g)$. Thus ${\cal
O}(s)(\nu _1)\ne {\cal O}(s)(\nu _2)$ (or ${\cal O}_h(\nu _1)\ne
{\cal O}_h(\nu _2)$, ${\cal O}_l(\nu _1)\ne {\cal O}_l(\nu _2)$,
${\cal O}_{l,h}(\nu _1)\ne {\cal O}_{l,h}(\nu _2)$ correspondingly).

\par {\bf 20. Corollary.} {\it The functors $I$, $I_h$, $I_l$ and $I_{l,h}$
are monomorphic.}
\par {\bf Proof.} This follows from Proposition 19.

\par {\bf 21. Proposition.} {\it The functors ${\cal O}$, ${\cal O}_h$, ${\cal O}_l$
and ${\cal O}_{l,h}$ are epimorphic.}
\par {\bf Proof.} Suppose that $f: X\to Y$ is a surjective mapping,
$\nu \in {\cal O}(Y,K)$ (or in ${\cal O}_h(Y,K)$, ${\cal O}_l(Y,K)$,
${\cal O}_{l,h}(Y,K)$ respectively). The set $L$ of all mappings
$g\circ f: X\to K$ with $g\in C(Y,K)$ (or in $C_+(Y,K)$
correspondingly) is the $\cal A$-subset or the left module over $K$
in $C(X,K)$ (or in $C_+(X,K)$). Put $\mu (g\circ f) = \nu (g)$. This
functional has an extension from $L$ to a functional $\mu \in {\cal
O}(X,K)$ (or in ${\cal O}_h(X,K)$, ${\cal O}_l(X,K)$, ${\cal
O}_{l,h}(X,K)$ correspondingly) due to Lemma 14.
\par {\bf 22. Lemma.} {\it Let $L$ be a submodule over $K$ of
$C(X,K)$ or $C_+(X,K)$ relative to the operations $\vee $, $\wedge
$, $\bigodot $ and containing all constant mappings $g^c: X\to K$,
where $c\in K$. Let also $\nu : L\to K$ be an idempotent (left
homogeneous) functional. For each $f\in C(X,K)\setminus L$ or
$C_+(X,K)\setminus L$ there exists an idempotent (left homogeneous)
extension $\mu _M$ of $\nu $ on a minimal submodule $M$ containing
$L$ and $f$.}
\par {\bf Proof.} For each $g\in M$ one puts
\par $(1)$ $\nu (g) = \inf \{ \nu (h): h\le g, h\in L \} $. \\
Therefore, $\nu (g_1)\le \nu (g_2)$ for each $g_1\le g_2\in M$. Then
\par $\nu (g^c\bigodot g) = \inf \{ \nu (h): h\in L,  g^c\bigodot
g\ge h \} = $\par $\inf \{ \nu (g^c\bigodot q): q\in L,  g^c\bigodot
g\ge g^c\bigodot q \}  = c\bigodot \inf \{ \nu (q): q\in L,  q\le g
\} = c\bigodot \nu (g)$ and \par $\nu (g\bigodot g^c) = \inf \{ \nu
(h): h\in L, g\bigodot g^c\ge h \} = \inf \{ \nu (q\bigodot g^c):
q\in L, q\bigodot g^c\le g\bigodot g^c \}$\par $ = \inf \{ \nu (q):
q\in L, q\le g \} \bigodot c =  \nu (g)\bigodot c$. \par On the
other hand for each $g_1, g_2\in M$ one gets
\par $\nu (g_1)\vee \nu (g_2) = \inf \{ \nu (g): g\in L, g_1\ge g \}
\vee \inf \{ \nu (q): q\in L, g_2\ge q \}$ \\ $ = \inf \{ \nu
(g)\vee \nu (q): g, q \in L, g_1\ge g, g_2\ge q \}  \ge \inf \{ \nu
(g\vee q): g, q \in L, g_1\vee g_2\ge g\vee q \} = \nu (g_1\vee
g_2)$.
\par From the inequalities $g_k\le g_1\vee g_2$ for $k=1$ and $k=2$ it follows,
that $\nu (g_k)\le \nu (g_1\vee g_2)$, consequently, $\nu (g_1)\vee
\nu (g_2)  = \nu (g_1\vee g_2)$. Then
\par $\nu (g_1)\wedge \nu (g_2) = \inf \{ \nu (g): g\in L, g_1\ge g \}
\wedge \inf \{ \nu (q): q\in L, g_2\ge q \}$ \\ $ = \inf \{ \nu
(g)\wedge \nu (q): g, q \in L, g_1\ge g, g_2\ge q \}  \le \inf \{
\nu (g\wedge q): g, q \in L, g_1\wedge g_2\ge g\wedge q \} = \nu
(g_1\wedge g_2)$.
\par But $\nu (g_k)\ge \nu (g_1\wedge g_2)$, since $g_k\ge g_1\wedge g_2$ for $k=1$ and $k=2$,
consequently, $\nu (g_1)\wedge \nu (g_2)  = \nu (g_1\wedge g_2)$. If
$\nu $ is left homogeneous, then $\inf \{ \nu (bh): bh \le bg, h\in
L \} = \inf \{ \nu (bh): h\le g, h\in L \} = b \inf \{ \nu (h): h\le
g, h \in L \} $ for each $b\in K$, consequently, $\nu $ is left
homogeoeneous on $M$.
\par {\bf 23. Lemma.} {\it If suppositions of Lemma 22 are
satisfied, then there exists an idempotent (left homogeneous)
functional $\lambda $ on $C(X,K)$ or $C_+(X,K)$ respectively such
that $\lambda |_L=\nu $.}
\par {\bf Proof.} The family of all extensions $(M,\mu _M)$ of $\nu $ on submodules
of $C(X,K)$ or $C_+(X,K)$ respectively is partially ordered by
inclusion: $(M,\mu _M)\le (N,\mu _N)$ if and only if $M\subset N$
and $\nu _N|_M=\nu _M$. In view of the Kuratowski-Zorn lemma
\cite{kunenb} there exists the maximal submodule $P$ in $C(X,K)$ or
$C_+(X,K)$ correspondingly and an idempotent extension $\nu _P$ of
$\nu $ on $P$. If $P\ne C(X,K)$ or $C_+(X,K)$ correspondingly by
Lemma 22 this functional $\nu _P$ could be extended on a module $L$
containing $P$ and some $g\in C(X,K)\setminus P$ or in
$C(X,K)_+\setminus P$ respectively. This contradicts the maximality
of $(P,\nu _P)$. Thus $P=C(X,K)$ or $C_+(X,K)$ correspondingly.
\par {\bf 24. Proposition.} {\it The functors $I$, $I_l$ and $I_h$, $I_{l,h}$
are epimorphic.}
\par {\bf Proof.} Let $f: X\to Y$ be epimorphic. We consider the set $L$
of all mappings $g\circ f: X\to K$ such that $g\in C(Y,K)$ or
$C_+(Y,K)$. Then $L$ is a submodule of $C(X,K)$ or $C_+(X,K)$
relative to the operations $\vee $, $\wedge $, $\bigodot $ and $L$
contains all constant mappings $g^c: X\to K$, where $c\in K$. Put
$\mu (g\circ f) = \nu (g)$ for $\nu \in I(X,K)$ or in $I_l(X,K)$,
$I_h(X,K)$ or $I_{l,h}(X,K)$. In view of Lemma 23 there is an
extension of $\mu $ from $L$ onto $C(Y,K)$ or $C_+(Y,K)$ such that
$\mu \in I(Y,K)$ or in $I_l(Y,K)$, $I_h(Y,K)$ or $I_{l,h}(Y,K)$
correspondingly.
\par {\bf 25. Definition.} It is said that $\nu \in {\cal O}(X,K)$ (or
$\nu \in {\cal O}_l(X,K)$) is supported on a subset $E$ in $X$, if
$\nu (f)=0$ for each $f\in C(X,K)$ or in $C_+(X,K)$ such that
$f|_E\equiv 0$. A support of $\nu $ is the intersection of all
subsets in $X$ on which $\nu $ is supported.
\par {\bf 26. Proposition.} {\it Let $\nu \in {\cal O}(X,K)$ or in
${\cal O}_l(X,K)$. Then $\nu $ is supported on $E\subset X$ if and
only if $\nu (f)=\nu (g)$ for each $f, g\in C(X,K)$ or in $C_+(X,K)$
correspondingly such that $f|_E \equiv g|_E$. Moreover, $E$ is a
support of $\nu $ if and only if $\nu $ is supported on $E$ and for
each proper subset $F$ in $E$, i.e. $F\subset E$ with $F\ne E$,
there are $f, h \in C(X,K)$ or in $C_+(X,K)$ respectively with $f|_F
\equiv h|_F$ such that $\nu (f)\ne \nu (h)$.}
\par {\bf Proof.} Consider $\nu \in {\cal O}(X,K)$ such that $\nu (f)=\nu (g)$
for each functions $f, g: X\to K$ with $f| _E=g|_E$. A functional
$\nu $ induces a functional $\lambda \in {\cal O}(E,K)$ such that
$\lambda (h)=\nu (h)$ for each $h\in C(X,K)$ with $h|_{X\setminus
E}=0$. Denote by $id$ the identity embedding of $E$ into $X$. Each
function $t: E\to K$ has an extension on $X$ with values in $K$.
Then ${\cal O}(id)(\lambda )=\nu $, since $\nu (g^0)=0$ and hence
$\nu (s)=0$ for each $s\in C(X,K)$ such that $s|_E\equiv 0$.
\par If $\nu \in {\cal O}(X,K)$ and $\nu $ is supported on $E$, then by
Definition 25 there exists a functional $\lambda \in {\cal O}(E,K)$
such that ${\cal O}(id)(\lambda )=\nu $. Therefore $\nu (f)=\lambda
(f|_E)=\lambda (g|_E)=\nu (g)$ for each functions $f, g \in C(X,K)$
such that $f| _E=g|_E$.
\par If $E$ is a support of $\nu $, then by the definition this
implies that $\nu $ is supported on $E$. Suppose that $F\subset E$,
$F\ne E$ and for each $f, g \in C(X,K)$ with $f|_F\equiv g|_F$ the
equality $\nu (f)=\nu (g)$ is satisfied, then a support of $\nu $ is
contained in $F$, hence $E$ is not a support of $\nu $. This is the
contradiction, hence there are $f, g \in C(X,K)$ with $f|_F\equiv
g|_F$ such that $\nu (f)\ne \nu (g)$.
\par If $\nu $ is supported on $E$ and for each
proper subset $F$ in $E$ there are $f, h \in C(X,K)$ with $f|_F
\equiv h|_F$ such that $\nu (f)\ne \nu (h)$, then $\nu $ is not
supported on any such proper subset $F$, consequently, each subset
$G$ in $X$ on which $\nu $ is supported contains $E$, i.e. $ ~
E\subset G$. Thus $E$ is the support of $\nu $. \par One can put
\par $(1)$ $t_1(x) = \sup \{ t(y): y\in E, y\le x \} $ for each
$x\in  X$. In the case ${\cal O}_l$ the proof is analogous, since
each nondecreasing mapping $t: E\to K$ has a nondecreasing extension
$t_1$ on $X$ with values in $K$, when $X\in Ob ({\cal S}_l)$,
$E\subset X$.
\par {\bf 27. Proposition.} {\it The functors ${\cal O}$, $I$, ${\cal O}_l$,
$I_l$, ${\cal O}_{l,h}$, $I_{l,h}$ preserve intersections.}
\par {\bf Proof.} If $E$ is a subset in $X$, then
there is the natural embedding $C(E,K)\hookrightarrow C(X,K)$ (or
$C_+(E,K)\hookrightarrow C_+(X,K)$, when $X\in Ob({\cal S}_l)$ due
to Formula 26$(1)$). Therefore, ${\cal O}(E\cap F,K)\subset {\cal
O}(E,K)\cap {\cal O}(F,K)$ (or ${\cal O}_l(E\cap F,K)\subset {\cal
O}_l(E,K)\cap {\cal O}_l(F,K)$ respectively). For any subsets $E$
and $F$ in $X$ and each functions $f, g \in C(X,K)$ (or $C_+(X,K)$)
with $f|_{E\cap F} \equiv g|_{E\cap F}$ there exists a function
$h\in C(X,K)$ (or $C_+(X,K)$) such that $h|_E=f$ and $h|_F=g$.
Therefore $\nu (f)=\nu (h)$ and $\nu (g)=\nu (h)$ for each $\nu \in
{\cal O}(E,K)\cap {\cal O}(F,K)$ (or in ${\cal O}_l(E,K)\cap {\cal
O}_l(F,K)$). In view of Proposition 26 the functors ${\cal O}$ and
${\cal O}_l$ preserve intersections. This implies that the functors
$I$, $I_l$, ${\cal O}_{l,h}$ and $I_{l,h}$ also have this property.
\par {\bf 28. Proposition.} {\it Let $ \{ X_b; p^b_a; V \} =: P $ be
an inverse system of sets $X_b$, where $V$ is a directed set,
$p^b_a: X_b\to X_a$ is a mapping for each $a\le b\in V$, $p_b:
X=\lim P \to X_b$ is a projection. Then the mappings
\par $(1)$ $s=({\cal O}(p_b): b\in V) : {\cal O}(X,K)\to {\cal
O}(P,K)$ and $s_h=({\cal O}_h(p_b): b\in V) : {\cal O}_h(X,K)\to
{\cal O}_h(P,K)$
\par $(2)$ $t=(I(p_b): b\in V) : I(X,K)\to I(P,K)$ and
$t_h=(I_h(p_b): b\in V) : I_h(X,K)\to I_h(P,K)$ \\ are bijective and
surjective algebraic homomorphisms. Moreover, if $X_b\in Ob({\cal
S}_l)$ and $p^b_a$ is order-preserving for each $a<b\in V$, then the
mappings \par $(3)$ $s_l=({\cal O}_l(p_b): b\in V) : {\cal
O}_l(X,K)\to {\cal O}_l(P,K)$ and $s_{l,h}=({\cal O}_{l,h}(p_b):
b\in V) : {\cal O}_{l,h}(X,K)\to {\cal O}_{l,h}(P,K)$
\par $(4)$ $t_l=(I_l(p_b): b\in V) : I_l(X,K)\to I_l(P,K)$ and
$t_{l,h}=(I_{l,h}(p_b): b\in V) : I_{l,h}(X,K)\to I_{l,h}(P,K)$ \\
also are bijective and surjective algebraic homomorphisms.}
\par {\bf Proof.} We consider the inverse system ${\cal O}(P) = ({\cal O}(X_a);
{\cal O}(p^a_b); V \} $ and its limit space $Y=\lim {\cal O}(P)$.
Then ${\cal O}(p^b_a) {\cal O}(p_b) = {\cal O}(p_a)$  for each $a\le
b\in V$, since $p^b_a\circ p_b=p_a$. Let $q: {\cal O}(X,K)\to Y$
denote the limit map of the inverse mapping system $q = \lim \{
{\cal O}(p_a); {\cal O}(p^a_b); V \} $.
\par A functional $\nu $ is in ${\cal O}(X,K)$ if and only if
${\cal O}(p_a)(\nu ) \in {\cal O}(X_a,K)$ for each $a\in V$, since
\par $(5)$ $f\in C(X,K)$ if and only if $f = \lim \{ f_b; p^b_a; V \} $ and
\par $(6)$ ${\cal O}(p_a)(\nu )(f_a) = \nu (f_a\circ p_a)= \nu _a(f_a)$,
where $\nu _a\in {\cal O}(X_a,K)$, $f_b\in C(X_b,K)$, $f_b=f_a\circ
p^b_a$ for each $a\le b\in V$, $p^b_b=id$, $f(x)=\{ f_a\circ p_a(x):
~ a \in V \} \in \theta (K)$ for each $x= \{ x_a: a \in V \} \in X$,
where $ \{ x_a: a \in V \} $ is a thread of $P$ such that $x_a\in
X_a$, $p^b_a(x_b)=x_a$ for each $a\le b \in V$, $\theta : K\to K^X$
is an order-preserving algebraic embedding, $\theta (K)$ is
isomorphic with $K$.
\par If $\nu , \lambda \in
{\cal O}(X,K)$ are two different functionals, this means that a
function $f\in C(X,K)$ exists such that $\nu _1(f)\ne \nu _2(f)$.
This is equivalent to the following: there exists $a\in V$ such that
$({\cal O}(p_a)(\nu ))(f) \ne ({\cal O}(p_a)(\lambda ))(f)$. Thus
the mappings $s$ and analogously $t$ are surjective and bijective.
\par On the other hand, \par $(7)$ $\nu _b(f_b\vee g_b)= \nu _b(f_b)\vee \nu _b(g_b)$
and \par $(8)$ $\nu _b(f_b\wedge g_b)= \nu _b(f_b)\wedge \nu
_b(g_b)$ for each $b\in V$ and each $\nu _b\in I(X_b,K)$ and every
$f_b, g_b\in C(X_b,K)$ such that either $f_b(x)< g_b(x)$ or
$f_b(x)=g_b(x)$ or $g_b(x)<f_b(x)$ for each $x\in X_b$, also \par
$(9)$ $\nu _b(g^c\bigodot f_b) = c\bigodot \nu _b(f_b)$ and \par
$(10)$ $\nu _b(f_b\bigodot g^c)=\nu _b(f_b) \bigodot c$ for each
$c\in K$ and $f_b\in C(X_b,K)$. Taking the inverse limit in
Equalities $(5-10)$ gives the corresponding equalities for $\nu \in
I(X,K)$, where $\nu = \lim \{ \nu _a; I(p^b_a); V \} $, hence $t$ is
the algebraic homomorphism.
\par Analogously $s$ preserves Properties $(9,10)$, that is
$\lambda = \lim \{ \lambda _a; {\cal O}(p^b_a); V \} $ is weakly
additive, where $\lambda _b\in {\cal O}(X_b,K)$ for each $b\in V$.
Suppose that $f\le g\in C(X,K)$, then $f_b\le g_b$ for each $b\in V$
due to $(5)$. From $\lambda _b(f_b)\le \lambda _b(g_b)$ for each
$b\in V$, the inverse limit decomposition $\lambda = \lim \{ \lambda
_b; {\cal O}(p^b_a); V \} $ and $(6)$ it follows that $\lambda $ is
order-preserving. \par If $X_b\in Ob ({\cal S}_l)$ for each $b\in
V$, then $X$ is linearly ordered: $x= \{ x_b: b\in V \} \le y = \{
y_b: b\in V \} $ if and only if $x_b\le y_b$ for each $b\in V$,
where $x, y \in X$ are threads of the inverse system $P$ such that
$p^b_a(x_b)=x_a$ for each $a\le b\in V$. Since $p^b_a$ is
order-preserving for each $a\le b \in V$ and each $f_b$ is
non-decreasing, then $f$ is nondecreasing and hence $f\in C_+(X,K)$
for each $f=\lim \{ f_b; p^b_a; V \} $, where $f_b\in C_+(X_b,K)$
and $f_b=f_a\circ p^b_a$ for each $a\le b \in V$ and $x\in X$,
$f(x)= \{ f_a\circ p_a(x): a \in V \} $.
\par Moreover, $\nu \in {\cal O}_h(X,K)$ is left homogeneous if and
only if $\theta (p_a)(\nu )$ is left homogeneous for each $b\in V$,
since $({\cal O}_h(p_a)(\nu ))(f_a) = \nu (f_a\circ p_a) = \nu
_a(f_a)$.
\par {\bf 29. Lemma.}  { \it There exist covariant functors ${\cal
O}_2$, $I_2$, and ${\cal O}_{l,2}$, $I_{l,2}$ and ${\cal O}_{h,2}$,
$I_{h,2}$ and ${\cal O}_{l,h,2}$, $I_{l,h,2}$ in the categories
${\cal K}_w$ and ${\cal K}$ and ${\cal K}_{w,l}$ and ${\cal K}_l$
respectively.}
\par {\bf Proof.} If $K_1, K_2, K_3 \in Ob({\cal K}_w)$, $ ~ u\in Mor (K_1 ,K_2)$, $ ~ v\in
Mor (K_2,K_3)$, $ ~ \nu \in I(X,K_1)$, then $(I_2(vu)(\nu ))(f)=
v\circ u \circ \nu (f_1) = [I_2(v)(I_2(u)(\nu ))](f)$ for each
$f_1\in C(X,K_1)$ such that $f(x)=v\circ u\circ f_1(x)$ for each
$x\in X$, where $X\in Ob ({\cal S})$. That is $I_2(vu) = I_2(v)
I_2(u)$. Evidently, $I_2(id)=1$.
\par If $f(x)\le g(x)$, then $u(f(x))\le u(g(x))$, where $x\in X$,
$f, g \in C(X,K_1)$. Therefore, if $f\vee g$ or $f\wedge g$ exists
in $C(X,K_1)$, then $u(f\vee g)= u(f)\vee u(g)$ or $u(f\wedge g)=
u(f)\wedge u(g)$ in $C(X,K_2)$ respectively. If $f, g\in C(X,K_1)$,
then $u(f(x)+g(x)) = u(f(x))+u(g(x))$ for each $x\in X$,
particularly for $f=g^c$ or $g=g^c$, where $c\in K_1$. Therefore,
$u(g^c\bigodot g)= g^{u(c)} \bigodot u(g)$ and $u(g\bigodot g^c)=
u(g) \bigodot g^{u(c)}$. To each $\nu _n\in {\cal O}(X,K_n)$ and
$u\in Mor (K_n,K_{n+1})$ there corresponds a functional $u\circ \nu
_n$ on $({\hat {id}},{\hat u})(C(X,K_n))$, $ ~ ({\hat {id}},{\hat
u})(C(X,K_n))\hookrightarrow C(X,K_{n+1})$ (see \S 15). If $u:
K_n\to K_{n+1}$ is not an epimorphism, the image $({\hat {id}},{\hat
u})(C(X,K_n))$ is a proper submodule over $u(K_n)$ in
$C(X,K_{n+1})$.
\par If $K_n, K_{n+1} \in Ob ({\cal K})$ and $X\in Ob ({\cal
S}_l)$, $u\in Mor (K_n,K_{n+1})$, then $u: C_+(X,K_n)\to
C_+(X,K_{n+1})$ is a homomorphism. If $K_n, K_{n+1}\in Ob ({\cal
K}_l)$ and $X\in Ob ({\cal S}_l)$ (or $K_n, K_{n+1}\in Ob ({\cal
K}_{w,l})$ and $X\in Ob ({\cal S})$) and $\nu \in {\cal O}_h(X,K_n)$
or in $I_h(X,K_n)$, $u\in Mor (K_n,K_{n+1})$, then $u\circ \nu \in
{\cal O}_h(X,K_{n+1})$ or in $I_h(X,K_{n+1})$ respectively.
\par This and the definitions above imply that ${\cal O}_2(u): {\cal
O}(X,K_1)\to {\cal O}(X,K_2)$, $I_2(u): I(X,K_1)\to I(X,K_2)$ and
${\cal O}_{l,2}(u)$, $I_{l,2}(u)$ and ${\cal O}_{h,2}(u)$,
$I_{h,2}(u)$ and ${\cal O}_{l,h,2}(u)$, $I_{l,h,2}(u)$ are the
homomorphisms. Thus ${\cal O}_2: {\cal K}_w\to {\cal OK}$ and ${\cal
O}_{l,2}: {\cal K}\to {\cal O}_l{\cal K}$, $I_2: {\cal K}_w\to {\cal
IK}$ and $I_{l,2}: {\cal K}\to {\cal I}_l{\cal K}$, ${\cal O}_{h,2}:
{\cal K}_{w,l}\to {\cal O}_h{\cal K}$, $I_{h,2}: {\cal K}_{w,l}\to
{\cal I}_h{\cal K}$, ${\cal O}_{l,h,2}: {\cal K}_l\to {\cal
O}_{l,h}{\cal K}$ and ${\cal I}_{l,h,2}: {\cal K}_l\to {\cal
I}_{l,h}{\cal K}$ are the covariant functors on the categories
${\cal K}_w$, ${\cal K}$, ${\cal K}_{w,l}$ and ${\cal K}_l$
correspondingly with values in the categories of (skew) idempotent
functionals, when a set $X\in Ob({\cal S})$ or in $Ob({\cal S}_l)$
correspondingly is marked.
\par {\bf 30. Proposition.} {\it The bi-functors $I$ on ${\cal S}\times {\cal
K}_w$, $I_l$ on ${\cal S}_l\times {\cal K}$, $I_h$ on ${\cal
S}\times {\cal K}_{w,l}$ and $I_{l,h}$ on ${\cal S}_l\times {\cal
K}_l$ preserve pre-images.}
\par {\bf Proof.} In view of
Proposition 24 and Lemma 29 $I$, $I_l$, $I_h$ and $I_{l,h}$ are the
covariant bi-functors, i.e. the functors in ${\cal S}$ or ${\cal
S}_l$ and the functors in ${\cal K}_w$ or ${\cal K}$ or ${\cal
K}_{w,l}$ or ${\cal K}_l$ correspondingly as well. For any functor
$F$ the inclusion $F(f^{-1}(B))\subset (F(f))^{-1}(F(B))$ is
satisfied.
\par Suppose the contrary that $I$ does not preserve pre-images.
This means that there exist $X, Y\in Ob({\cal S})$ and $K_1, K_2 \in
Ob({\cal K}_w)$ or $X, Y\in Ob({\cal S}_l)$ and $K_1, K_2 \in
Ob({\cal K})$, $ ~ f\in Mor (X,Y)$, $ ~ u\in Mor (K_1,K_2)$,
$~A\subset X$ and $B\subset Y$, $~ A=F^{-1}(B)$, $~ \nu \in
I(X,K_1)$ such that $I(f,u)(\nu ) \in I(B,K_2)$ but $\nu \notin
I(f^{-1}(B),u^{-1}(K_2))$ (or $\nu \in I_l(X,K_1)$, $I_l(f,u)(\nu
)\in I_l(B,K_2)$ and $\nu \notin I_l(f^{-1}(B),u^{-1}(K_2))$
respectively). One can choose two functions $g, h \in C(X,K_1)$ such
that
\par $(1)$ $g|_A = h|_A$,
\par $(2)$ $0< c_1= u[\inf_{x\in X} g(x)]$, $0< c_2=u[\inf_{x\in X}
h(x)]$ and
\par $(3)$ $u[\nu (g)]\ne u[\nu (h)]$. \par There exist functions $s, t \in
C(X,K_1)$ such that \par $(4)$ $s|_A=g|_A$ and $t|_A=h|_A$, while
\par $(5)$ $s|_{X\setminus A}=t|_{X\setminus A}$ and \par $(6)$ $s(x)\le g(x)$ and
$s(x)\le h(x)$ for each $x\in X\setminus A$, where $g, h$ satisfy
Conditions $(1-3)$. There are also functions $q, r \in C(X,K_1)$
such that
\par $(7)$ $q|_{X\setminus A}=g|_{X\setminus A}$ and $r|_{X\setminus
A}=h|_{X\setminus A}$ with
\par $(8)$ $q(x)=r(x)$ and $q(x)\le c$ for each $x\in A$, where \par $(9)$ $c\in K_1$,
$c< \inf_{x\in X} g(x)$, $c<\inf _{x\in X} h(x)$ such that
$u(c)<c_1$ and $u(c)< c_2$. \par Evidently, $c_1\le u[\nu (g)]$ and
$c_2\le u[\nu (h)]$. Then \par $(10)$ $\nu (g) = \nu (s\vee q)= \nu
(s)\vee \nu (q)$ and \par $(11)$ $\nu (h) = \nu (t\vee r) =\nu
(t)\vee \nu (r)$ and $u[\nu (q)]\ne u[\nu (r)]$.
\par  On
the other hand, there are functions $q_1, r_1 \in C(Y,K_2)$, $ ~
q_2, r_2\in C(Y,K_1)$ such that $q_2\circ f=q$, $~ r_2\circ f=r$, $
~ u\circ q_2= q_1$, $ ~ u\circ r_2= r_1$ and $q_2|_B=r_2|_B$.
Therefore, from Properties $(7-10)$ it follows that \par $(12)$
$(I(f,u)(\nu ))(q_1) = u[\nu (q)]\le u(c)$ and $(I(f,u)(\nu ))(r_1)
= u[\nu (r)]\le u(c)$. The condition $s=t$ on $A$ and on $X\setminus
A$ imply that \par $(13)$ $\nu (s) = \nu (t)$. Therefore, \par
$(14)$ $u(\nu (g)) = u(\nu (s))\vee u(\nu (q))$ and $u(\nu (h)) =
u(\nu (t))\vee u(\nu (r))$, which follows from $(10,11)$. But
Formulas $(4-6,12-14)$ contradict the inequality $u[\nu (g)]\ne
u[\nu (h)]$, since $u$ is the order-preserving algebraic
homomorphism. Thus the bi-functors $I$ and $I_l$ preserve
pre-images. The proof in other cases is analogous.

\par {\bf 31. Corollary.} {\it If $\nu \in I(X,K)$ or $\nu \in I_l(X,K)$, $f\in Mor
(X,Y)$, $u\in Mor(K_1,K_2)$, where $X, Y\in Ob({\cal S})$ and $K_1,
K_2\in Ob({\cal K}_w)$ or $X, Y\in Ob({\cal S}_l)$ and $K_1, K_2\in
Ob({\cal K})$, then $supp (I(f,u)(\nu )) = f(supp (u[\nu ]))$ or
\par $supp (I_l(f,u)(\nu )) = f(supp (u[\nu ]))$ correspondingly.}

\par {\bf 32. Definitions.} Suppose that $Q$ is a category and $F,
G$ are two functors in $Q$. Suppose also that a transformation $p:
F\to G$ is defined for each $X\in Q$, that is a mapping $p_X:
F(X)\to G(X)$ is given. If $p_Y\circ F(f)=G(f)\circ p_X$ for each
mapping $f\in Mor(X,Y)$ and every objects $X, Y\in Ob(Q)$, then the
transformation $p = \{ p_X: X\in Ob (Q) \} $ is called natural.
\par If $T: Q\to Q$ is an endofunctor in a category $Q$ and there
are natural transformations the identity $\eta : 1_Q\to T$ and the
multiplication $\psi: T^2\to T$ satisfying the relations $\psi \circ
T\eta = \psi \circ \eta T = 1_T$ and $\psi\circ \psi T= \psi \circ
T\psi $, then one says that the triple ${\bf T} := (T, \eta , \psi)$
is a monad.

\par {\bf 33. Theorem.} {\it There are monads in the categories
${\cal S}\times {\cal K}_w$, ${\cal S}_l\times {\cal K}$, ${\cal
S}\times {\cal K}_{w,l}$ and ${\cal S}_l\times {\cal K}_l$.}
\par {\bf Proof.} Let ${\bar g} (\nu ) := \nu (g)$ for $g\in C(X,K)$
and $\nu \in I(X,K)$, where $X\in Ob ({\cal S})$ and $K\in Ob ({\cal
K}_w)$. Therefore, ${\bar g} : I(X,K)\to K$. Then
$$\overline{g^b\bigodot g} (\nu ) = \nu (g^b\bigodot g)= b\bigodot
\nu (g)= b\bigodot {\bar g}(\nu )\mbox{ and}$$
$$\overline{g\bigodot g^b} (\nu ) = \nu (g\bigodot g^b)= \nu (g)\bigodot
b= \overline{g}(\nu )\bigodot b,$$ where $g^b(x)=b$ for each $x\in
X$, that is
\par $(1)$ $\overline{g\bigodot g^b} ={\bar g}\bigodot
\overline{g^b}$ and $\overline{g^b\bigodot g} = \overline{g^b}
\bigodot {\bar g}$ \\ for each $g\in C(X,K)$ and $b\in K$. \par Then
we get $\overline{g\vee h} (\nu ) = \nu (g\vee h) = \nu (g)\vee \nu
(h) ={\bar g}(\nu )\vee {\bar h}(\nu )= ({\bar g}\vee {\bar h})(\nu
)$. Moreover, we deduce that $\overline{g\wedge h} (\nu ) = \nu
(g\wedge h) = \nu (g)\wedge \nu (h) = {\bar g}(\nu )\wedge {\bar
h}(\nu )= ({\bar g}\wedge {\bar h})(\nu )$. Thus \par $(2)$
$\overline{g\vee h} = {\bar g}\vee {\bar h}$ and $\overline{g\wedge
h} = {\bar g}\wedge {\bar h}$. \par If additionally $\nu $ is left
homogeneous and $K\in Ob ({\cal K}_{w,l})$, then $\overline{bg} =
\nu (bg) = b \nu (g) = b {\bar g}(\nu )$, hence $\overline{bg} =
b\bar g$ for every $b\in K$ and $g\in C(X,K)$.
\par For $\lambda \in I(I(X,K),K)$ we put $\xi _{X,K}(\lambda
)(g)= \lambda ({\bar g})$ for each $g\in C(X,K)$. Then $\xi
_{X,K}(\lambda )(g^b)= \lambda (\overline{g^b}) = \lambda (q^b) =b$,
where $q^b: I(X,K)\to K$ denotes the constant mapping $q^b(y)=b$ for
each $y\in I(X,K)$. From Formulas $(1)$ it follows that $$\xi
_{X,K}(\lambda )(g^b\bigodot g) = \lambda (\overline{g^b\bigodot g})
= \lambda (b\bigodot {\bar g}) = b \bigodot \lambda ({\bar g}) = b
\bigodot \xi _{X,K}(\lambda )(g)\mbox{  and}$$
$$\xi _{X,K}(\lambda )(g\bigodot g^b) =
\lambda (\overline{g\bigodot g^b}) = \lambda ({\bar g}\bigodot
\overline{g^b}) = \lambda ({\bar g}) \bigodot b = \xi _{X,K}(\lambda
)(g) \bigodot b.$$ On the other hand, from Formulas $(2)$ we get
that
$$\xi _{X,K}(\nu )(g\vee h) = \nu (\overline{g\vee h})= \nu ({\bar
g}\vee {\bar h}) = \nu ({\bar g})\vee \nu ({\bar h}) = \xi
_{X,K}(\nu )(g)\vee \xi _{X,K}(\nu )(h)\mbox{  and}$$
$$\xi _{X,K}(\nu )(g\wedge h) = \nu (\overline{g\wedge h})= \nu ({\bar g}\wedge
{\bar h}) = \nu ({\bar g})\wedge \nu ({\bar h}) = \xi _{X,K}(\nu
)(g)\wedge \xi _{X,K}(\nu )(h)$$ for each $b\in K$, $g, h \in
C(X,K)$. Thus $\xi _{X,K}: I(I(X,K),K)\to I(X,K)$.
\par If $\lambda \in I_h(I_h(X,K),K)$ for $K\in Ob ({\cal
K}_{w,l})$, then $\xi _{X,K}(\lambda )(bg) = \lambda (\overline{bg})
= \lambda (b {\bar g}) = b \lambda ({\bar g})$, hence $\xi _{X,K}:
I_h(I_h(X,K),K)\to I_h(X,K)$.
 Analogously is defined the mapping $\xi _{X,K}: {\cal O}({\cal O}(X,K),K) \to {\cal
O}(X,K)$ for each $X\in Ob({\cal S})$ and $K\in {\cal K}_w$, also
$\xi _{X,K}: {\cal O}_l({\cal O}_l(X,K),K) \to {\cal O}_l(X,K)$,
$\xi _{X,K}: I_l(I_l(X,K),K)\to I_l(X,K)$ for each $X\in Ob({\cal
S}_l)$ and $K\in {\cal K}$, $\xi _{X,K}: I_h(I_h(X,K),K)\to
I_h(X,K)$ for $X\in Ob ({\cal S})$ and $K\in {\cal K}_{w,l}$, $\xi
_{X,K}: I_{l,h}(I_{l,h}(X,K),K)\to I_{l,h}(X,K)$ for $X\in Ob ({\cal
S}_l)$ and $K\in {\cal K}_l$. One also puts $\eta : Id_Q\to {\cal
O}$ or $\eta : Id_Q\to I$ for $Q={\cal S}\times {\cal K}_w$, also
$\eta : Id_Q\to {\cal O}_l$ or $\eta : Id_Q\to I_l$ for $Q = {\cal
S}_l\times {\cal K}$ correspondingly.
\par Next we verify that the transformations $\eta $ and $\xi $ are natural for each
$f\in Mor(X\times K_1,Y\times K_2)$, i.e. $f=(s,u)$, $s\in
Mor(X,Y)$, $u\in Mor (K_1,K_2)$: $$\eta _{(Y,K_2)}\circ {\cal
O}((s,u)) = {\cal O}(id_Y,id_{K_2})) \circ {\cal O}((s,u))$$ $$ =
{\cal O}((s,u)) = {\cal O}((s,u))\circ {\cal O}(id_X,id_{K_1})) =
{\cal O}((s,u))\circ \eta _{(X,K_1)},$$
$$\xi _{(Y,K_2)}\circ {\cal
O}((s,u))[{\cal O}^2(X,K_1)] = \xi _{(Y,K_2)}({\cal O}({\bar
s},{\bar u})[{\cal O}(X,K_1)])$$ $$ = {\cal O}((s,u))\circ \eta
_{(X,K_1)}[{\cal O}^2(X,K_1)]),$$ where ${\cal O}^{m+1}(X,K) :=
{\cal O}({\cal O}^m(X,K),K)$ for each natural number $m$ (see also
\S 15 and Proposition 30).
\par For each $\nu \in {\cal O}(X,K)$ and $g\in C(X,K)$
one gets $$\xi _{X,K} \circ \eta _{({\cal O}(X,K),K)} (\nu ) (g) =
\eta _{({\cal O}(X,K),K)} (\nu )({\bar g}) = {\bar g}(\nu )=\nu (g)
\mbox{ and}$$
$$\xi _{X,K} \circ {\cal O}(\eta _{(X,K)}) (\nu ) (g) =
({\cal O}(\eta _{(X,K)} (\nu ))({\bar g}) = \nu ({\bar g}\circ \eta
_{(X,K)}) =\nu (g).$$ Let now $\tau \in {\cal O}^3(X,K)$ and $g\in
C(X,K)$, then $$\xi _{(X,K)} \circ \xi _{{\cal O}(X,K)}(\tau )(g) =
(\xi_{{\cal O}(X,K)} (\tau )) ({\bar g}) = \tau ({\bar {\bar g}})
\mbox{ and}$$  $$\xi _{(X,K)} \circ {\cal O}(\xi _{(X,K)} )(\tau
)(g) = ({\cal O}(\xi _{(X,K)} )(\tau )) ({\bar g}) = \tau ({\bar g}
\circ \xi _{(X,K)}) = \tau ({\bar {\bar g}}),$$ where ${\bar {\bar
g}} \in C({\cal O}^2(X,K),K)$ is prescribed by the formula $({\bar
{\bar g}} )(\nu ) = \nu ({\bar g})$ for each $\nu \in {\cal
O}^2(X,K)$. Thus ${\bf O} := ({\cal O},\eta , \xi )$ is the monad.
Since $I$ is the restriction of the functor ${\cal O}$, the triple
${\bf I} := (I, \eta , \xi )$ is the monad in the category ${\cal
S}\times {\cal K}_w$ as well. Analogously ${\bf O}_l := ({\cal
O}_l,\eta , \xi )$ and ${\bf I}_l := (I_l, \eta , \xi )$ form the
monads in the category ${\cal S}_l\times {\cal K}$; ${\bf O}_h=
({\cal O}_h, \eta ,\xi )$ and ${\bf I}_h=(I_h,\eta , \xi )$ are the
monads in ${\cal S}\times {\cal K}_{w,l}$; ${\bf O}_{l,h} = ({\cal
O}_{l,h}, \eta ,\xi )$ and ${\bf I}_{l,h}=(I_{l,h},\eta , \xi )$ are
the monads in ${\cal S}_l\times {\cal K}_l$.
\par {\bf 34. Proposition.} {\it If a sequence \par $(1)$ $...\to K_n\to
K_{n+1}\to K_{n+2}\to ... $ in ${\cal K}_w$  (or in ${\cal K}$) is
exact, then sequences \par $(2)$ $...\to {\cal O}_2(X,K_n)\to {\cal
O}_2(X,K_{n+1})\to {\cal O}_2(X,K_{n+2})\to ... $ and \par $(3)$
$...\to I_2(X,K_n)\to I_2(X,K_{n+1})\to I_2(X,K_{n+2})\to ... $
\\ are exact (analogously for ${\cal O}_{l,2}$ and $I_{l,2}$
correspondingly).}
\par {\bf Proof.} A sequence \par $...\to K_n\to
K_{n+1}\to K_{n+2}\to ... $ is exact means that $s_n(K_n)= ker
(s_{n+1})$ for each $n$, where $s_n: K_n\to K_{n+1}$ is an
order-preserving algebraic homomorphism, $ker (s_{n+1}) =
s_{n+1}^{-1}(0)$. Each homomorphism $s_n$ induces the homomorphism
${\bf s}_n: C(X,K_n)\to C(X,K_{n+1})$ point-wise $({\bf
s}_n(f))(x)=s_n(f(x))$ for each $x\in X$. Therefore, ${\bf
s}_n(f\vee g)={\bf s}_n(f)\vee {\bf s}_n(g)$ or ${\bf s}_n(f\wedge
g)={\bf s}_n(f)\wedge {\bf s}_n(g)$, when $f\vee g$ or $f\wedge g$
exists, where $f, g \in C(X,K_n)$. Moreover, $({\bf s}_n(f+g))(x) =
{\bf s}_n(f(x)+g(x))=s_n(f(x))+s_n(g(x))=[{\bf s}_n(f)+{\bf
s}_n(g)](x)$ and $[{\bf s}_n(fg)](x)=
s_n(f(x)g(x))=s_n(f(x))s_n(g(x))=[({\bf s}_n(f))({\bf s}_n(g))](x)$,
consequently, ${\bf s}_n(C(X,K_n))={\bf s}_{n+1}^{-1}(0)$, since
$f_{n+2}\in C(X,K_{n+2})$ is zero if and only if $f_{n+2}(x)=0$ for
each $x\in X$. Thus the sequence
\par $...\to C(X,K_n)\to C(X,K_{n+1})\to C(X,K_{n+2})\to ... $ is
exact. \par Then a functional $\lambda _{n+2}\in {\cal
O}(X,K_{n+2})$ is zero on ${\bf s}_{n+1}(C(X,K_{n+1}))$ if and only
if $\lambda _{n+2}(f_{n+2})=0$ for each $f_{n+2}\in {\bf
s}_{n+1}(C(X,K_{n+1}))$. Therefore, ${\bf s}_{n+1}(\lambda
_{n+1})=0=\lambda _{n+2}$ on ${\bf s}_{n+1}[{\bf s}_n(C(X,K_n))]$ if
and only if $\lambda _{n+1}(f_{n+1})\in s_n(K_n)$ for each
$f_{n+1}\in {\bf s}_n(C(X,K_n))$. But ${\bf s}_{n+1}[{\bf
s}_n(C(X,K_n))]\subset {\bf s}_{n+1}(C(X,K_{n+1}))$, consequently,
${\cal O}_2(s_n)=ker {\cal O}_2(s_{n+1})$. Thus the sequences
$(2,3)$ are exact, analogously for other functors $I_2$, ${\cal
O}_{l,2}$ and $I_{2,l}$.

\par {\bf 35. Lemma.} {\it Let $G$ be a groupoid with the
unit acting on a set $X$ such that to each element $g\in G$ a
mapping $v_g: X\to X$ corresponds having the properties
\par $(1)$ $v_gv_h=v_{gh}$ for each $g, h \in G$ and \par $(2)$ $v_e=id$, where $e\in
G$ is the unit element, $id(x)=x$ for each $x\in X$. If $K$ is a
quasiring with the associative sub-quasiring $L$, $ ~ L\supset \{ 0,
1 \} $, such that \par $(3)$ $a(bc)=(ab)c$ for each $a, b\in L$ and
$c\in K$, a mapping $\rho : G^2\to L\setminus \{ 0 \} $ satisfies
the cocycle condition
\par $(4)$ $\rho (g,x)\rho (h,v_gx) = \rho (gh,x)$ and \par $(5)$
$\rho (e,x)=1\in K$ for each $g, h\in G$ and $x\in X$, then
\par $(6)$ $T_g f(x) := \rho (g,x) {\hat v}_gf(x)$ is a representation of
$G$ by mappings $T_g$ of $C(X,K)$ into $C(X,K)$, where $f\in
C(X,K)$, $ ~ {\hat v}_gf(x) := f(v_g(x))$ for each $g\in G$ and
$x\in X$.}
\par {\bf Proof.} For each $g, h \in G$ one has $T_g(T_hf(x))= \rho (g,x)
{\hat v}_g [\rho (h,x) {\hat v}_hf(x)] = \rho (gh,x){\hat
v}_{gh}f(x) = T_{gh} f(x)$, hence $T_gT_h=T_{gh}$. Moreover,
$T_ef=f$, since $v_e=id$ and $\rho (e,x)=1$, i.e. $T_e=I$ is the
unit operator on $C(X,K)$.
\par The mappings $T_g$ are (may be) generally non-linear relative to $K$.
If $K$ is commutative, distributive and associative, then $T_g$ are
$K$-linear on $C(X,K)$.
\par {\bf 36. Definition.} A
functional $\nu $ on $C(X,K)$ or $C_+(X,K)$ we call semi-idempotent,
if it satisfies the property:
\par $(1)$ $\nu (g+ f) = \nu (g)+ \nu (f)$ for each $f, g \in C(X,K)$
or $C_+(X,K)$ respectively, where $(g+f)(x) = g(x)+ f(x)$ for each
$x\in X$.
\par Suppose that $G$ is a groupoid with the
unit acting on a set $X$ and satisfying Conditions 35$(1,2)$. A
functional $\lambda $ on $C(X,K)$ or $C_+(X,K)$ we call
$(T,G)$-invariant if \par $(2)$ ${\hat T}_g\lambda = \lambda $,
where $({\hat T}_g\lambda )(f) := \lambda (T_gf)$ for each $g\in G$
and $f$ in $C(X,K)$ or $C_+(X,K)$ correspondingly. \par Let
$S_+(G,K)$ denote the family of all semi-idempotent functionals,
when $K$ is commutative and associative relative to the addition for
$(G,K)$, let also $S_{\vee }(G,K)$ (or $S_{\wedge }(G,K)$) denote
the family of all functionals satisfying Conditions 8$(4)$ (or
8$(5)$ correspondingly) for general $K$. Denote by $H_+(G,K)$ (or
$H_{\vee }(G,K)$ or $H_{\wedge }(G,K)$) the family of all
$G$-invariant semi-idempotent (or in $S_{\vee }(G,K)$ or in
$S_{\wedge }(G,K)$ correspondingly) functionals for $(X,K)$, when
$X=G$ as a set. We supply these families with the operations of the
addition
\par $(3)$ $\nu (f)+_i \lambda (f)=: (\nu +_i  \lambda )(f)$ in $S_j(G,K)$
for $i=1, 2, 3$ and $j= +, \vee , \wedge $ respectively and the
multiplication being the convolution of functionals
\par $(4)$ $\nu *\lambda (f)=\nu (\lambda (T_gf))$ in $S_j(G,K)$,
where $g\in G$, $~j\in \{ +, \vee , \wedge \} $. \par Then we put
$H_h(G,K)$, $S_h(G,K)$, $H_{\vee ,h}(G,K)$, $S_{\vee ,h}(G,K)$,
$H_{\wedge ,h}(G,K)$ and $S_{\wedge ,h}(G,K)$ for the subsets of all
left homogeneous functionals in $H_+(G,K)$, $S_+(G,K)$, $H_{\vee
}(G,K)$, $S_{\vee }(G,K)$, $H_{\wedge }(G,K)$, $S_{\wedge }(G,K)$
correspondingly.
\par {\bf 37. Proposition.} {\it If $\nu $ is a $(T,G)$-invariant
semi-idempotent functional, then its support is contained in
$\bigcap _{n=1}^{\infty } T^n(X) $, where $$T(A):= \bigcup_{g\in G}
supp (\rho (g,x) {\hat v}_g(\chi _A(x)))$$ for a subset $A$ in $X$.
Moreover, if $K$ has not divisors of zero a support of $\nu $ is
$G$-invariant and contained in $\bigcap _{n=1}^{\infty } P^n(X)$,
where $$P(X)= \bigcup_{g\in G} v_g(X).$$}
\par {\bf Proof.} If $\nu (f)\ne 0$, then $\nu (T_gf)\ne 0$ for each
$g\in G$, when a functional $\nu $ is $(T,G)$-invariant. On the
other hand, if $supp (f)\subset supp (\nu )$, then $supp (\rho (g,
x){\hat v}_gf(x)) \subset supp (\nu )$. At the same time,
$\bigcup_{g\in G} supp (T_gf) \subset \bigcup_{g\in G} supp ({\hat
v}_gf)$, since $\rho (g,x)\in L\setminus \{ 0 \} $ for each $g\in G$
and $x\in X$. Taking $f=\chi _{supp (\nu )}$ we get $supp (\nu
)\subset T(supp (\nu ))\subset T(X)$, hence by induction $supp (\nu
)\subset T^n(X)$ for each natural number $n$.
\par If $K$ has not divisors of zero, then $supp ({\hat T}_g\nu ) = {\hat v}_g supp(\nu )
\subset supp (\nu )$ for each $g\in G$, hence $\bigcup_{g\in G}
{\hat v}_g supp(\nu ) = supp (\nu )$, since $e\in G$ and $\nu
_e=id$. That is $supp (\nu )$ is $G$-invariant. Since $supp (\nu
)\subset X$, then $supp (\nu ) \subset P(X)$ and by induction $supp
(\nu )\subset P^n(X)$ for each natural number $n$.
\par {\bf 38. Proposition.} {\it If $G$ is a groupoid with a unit or
a monoid, then $S_+(G,K)$, $S_{\vee }(G,K)$ and $S_{\wedge }(G,K)$
for general $T_g$ and $K$ (or $S_h(G,K)$, $S_{h,\vee }(G,K)$ and
$S_{h, \wedge }(G,K)$ for $T_g\equiv {\hat v}_g$ or when $K$ is
commutative and associative relative to the multiplication) supplied
with the convolution 36$(4)$ as the multiplication operation are
groupoids with a unit or monoids correspondingly.}
\par {\bf Proof.} Certainly, the definitions above imply the inclusion
$S_h(G,K)\subset S_+(G,K)$. If $\nu , \lambda \in I_h(G,K)$, then
$(\nu *\lambda )(bf) = \nu (\lambda (T_g(bf))) = \nu (b\lambda
(T_gf)) = b ((\nu *\lambda )(f))$, when $T_g\equiv {\hat v}_g$ or
$K$ is commutative and associative relative to the multiplication.
We mention that the Dirac functional $\delta _e$ belongs to
$S_h(G,K)$ and has the property $\nu *\delta _e = \delta _e*\nu =\nu
$ for each $\nu \in S(G,K)$, where $e$ is a unit element in $G$.
Thus $\delta _e$ is the neutral element in $S(G,K)$.
\par For a monoid $G$ one has ${\hat v}_s({\hat v}_uf(x)) = f(s(ux))
= f((su)x) = {\hat v}_{su}f(x)$ for each $f\in C(G,K)$ and $s, u,
x\in G$.
\par If $G$ is
a monoid, then $(\nu *(\lambda * \phi )) (f) = \nu ^u ((\lambda *
\phi )(T_uf)) = \nu ^u (\lambda ^s(\phi (T_sT_uf)) = \nu ^u (\lambda
^s(\phi (T_{su}f))) = (\nu * \lambda )^{su} (\phi (T_{su}f)) = [(\nu
*\lambda )*\phi ](f)$ for every $f\in C(G,K)$ and $u, s\in G$ and
$\nu , \lambda , \phi \in S_j(G,K)$, where $\nu ^u(h)$ means that a
functional $\nu $ on a function $h$ acts by the variable $u\in G$,
consequently, $\nu
*(\lambda * \phi ) = (\nu *\lambda )* \phi $. Thus the family $S_j(G,K)$ is
associative, when $G$ is associative, where $j \in \{ +, \vee
,\wedge , h, (h,\vee ), (h,\wedge ) \} $ for the corresponding $T_g$
and $K$.

\par {\bf 39. Theorem.} {\it If $G$ is a groupoid with a unit or
a monoid, then $S_+(G,K)$ (for $K$ commutative and associative
relative to $+$), $S_{\vee }(G,K)$ and $S_{\wedge }(G,K)$ for
general $T_g$ (or $S_{\vee ,h}(G,K)$ and $S_{\wedge ,h}(G,K)$ for
$T_g \equiv {\hat v}_g$ or when $K$ is commutative and associative
relative to the multiplication) are quasirings or semirings
correspondingly.}
\par {\bf Proof.}  If $f, g \in C(X,K)$ or in $C_+(X,K)$ and
$f\vee g$ or $f\wedge g$ exists (see Condition $(3)$ in Lemma 6),
$\nu , \lambda $ are functionals satisfying Condition 8$(4)$ or
8$(5)$ respectively, then \par $(1)$ $(\nu +_i \lambda )(f+_i g) =
\nu (f+_i g) +_i \lambda (f+_i g) = (\nu (f)+_i \nu (g))+_i (\lambda
(f)+_i \lambda (g)) = (\nu (f)+_i \lambda (f)) +_i (\nu (g)+_i
\lambda (g)) = (\nu +_i \lambda )(f) +_i (\nu +_i \lambda )(g)$ \\
for $i=1, 2, 3$, where $+_1=+$, $+_2=\vee $, $+_3=\wedge $. That is,
the functional $\nu +_i \lambda $ satisfies Property 36$(1)$ for
$i=1$ or 8$(4)$ for $i=2$ or 8$(5)$ when $i=3$ correspondingly. If
additionally $\nu $ and $\lambda $ are left homogeneous, then \par
$(2)$ $(\nu +_i \lambda )(bf) = \nu (bf)+_i \lambda (bf) = b\nu
(f)+_i b \lambda (f) = b (\nu +_i \lambda )(f)$ \\ for each $b\in
K$.
\par On the other hand, we deduce that \par $((\nu _1+_i\nu
_2)*\lambda )(f)= (\nu _1+_i\nu _2) (\lambda (T_gf)) = \nu
_1(\lambda (T_gf)) +_i \nu _2(\lambda (T_gf))$\par $ = (\nu
_1*\lambda )(f) +_i (\nu _2*\lambda )(f)$ and
\par $(\lambda  *(\nu _1+_i\nu _2))(f)= \lambda ((\nu _1+_i\nu _2)
(T_gf)) = \lambda (\nu _1(T_gf)) +_i \lambda (\nu _2(T_gf))$\par $ =
(\lambda  *\nu _1)(f) +_i (\lambda  *\nu _2)(f)$ \\ for each $\nu
_1, \nu _2, \lambda \in S_j(G,K)$ and $f\in C(G,K)$ or in $C_+(G,K)$
correspondingly, for $i=1, 2, 3$ and $i=i(j)$ respectively, where
$+_1=+$, $+_2=\vee $ and $+_3=\wedge $. Thus, the right and left
distributive rules are satisfied: \par $(3)$ $(\nu _1+_i\nu _2)
* \lambda = \nu _1* \lambda +_i \nu _2* \lambda $ and
\par $(4)$ $\lambda  * (\nu _1+_i\nu
_2) = \lambda  * \nu _1 +_i \lambda  * \nu _2$ \\ for $i=1, 2, 3$
respectively.
\par Therefore, Formulas $(1-4)$ and Proposition 38 imply that
$S_+(G,K)$, $S_{\vee }(G,K)$, $S_{\wedge }(G,K)$, $S_{\vee ,h}(G,K)$
and $S_{\wedge ,h}(G,K)$ are left and right distributive quasirings
or semirings correspondingly.

\par {\bf 40. Theorem.} {\it If $G$ is a groupoid with a unit, $X=G$ as a set (see \S 36),
then $H_j(G,K)$ is an ideal in $S_j(G,K)$, where $j=+$ (for $K$
commutative and associative relative to $+$) or $j=\vee $ or
$j=\wedge $ or $j=(\vee ,h)$ or $j=(\wedge ,h)$ with $\rho
(u,x)\equiv 1$; $j=(\vee , h)$ or $j=(\wedge , h)$ for commutative
and associative $K$ relative to the multiplication with general
$T_u$.}
\par {\bf Proof.}  We mention that ${\hat T}_g(b_1\lambda
_1+_ib_2\lambda _2) (f) = b_1\lambda _1(T_gf) +_i b_2\lambda
_2(T_gf)$, where the operation denoted by the addition $+_i$ is
either $+$ or $\vee $ or $\wedge $ for $i=1$ or $i=2$ or $i=3$
correspondingly (and also below in this section), consequently,
$b_1\lambda _1+_ib_2\lambda _2\in H_j(G,K)$ for each $\lambda _1,
\lambda _2\in H_j(G,K)$ and $b_1, b_2\in K$, $i=i(j)$.
\par In Formula 36$(4)$ after the action of a functional $\lambda $
on a function $T_gf(x)$ of the variable $x$ one gets that $\lambda
(T_gf)=: h(g)$ is a function in the variable $g$ and $\nu $ is
acting on this function, i.e. $\nu
*\lambda (f)=\nu (h(x))$, where $x, g \in G$.
This implies that \par $\nu *(\lambda (f+_it)) = \nu * (\lambda (f)
+_i \lambda (t)) =  \nu (\lambda (T_gf) +_i \lambda (T_gt))$\par $ =
\nu ( \lambda (T_gf) )+_i \nu (\lambda (T_gt))= (\nu
*\lambda )(f)+_i (\nu *\lambda )(t)$ for $i=1, 2, 3$, \\ consequently, the
convolution operation maps from $S_j(G,K)^2$ into $S_j(G,K)$.
\par If $\lambda \in
H_j(G,K)$ and $\nu \in S_j(G,K)$, then \par $({\hat T}_s(\nu
* \lambda )) (f) = {\hat T}_s(\nu ^u(\lambda ^x (T_uf(x)))) = \nu ^u (\lambda ^x(T_s(T_uf(x)))$\par $ =
\nu ^u (\lambda ^x(T_uf(x)))) = (\nu *\lambda ) (f)$ and \par
$({\hat T}_s(\lambda *\nu )) (f) = {\hat T}_s(\lambda ^u(\nu
^x(T_uf(x))))= \lambda ^u
(\nu ^x(T_s(T_uf(x))))$\par $ = \lambda ^u(T_s(\nu ^x(T_uf(x)))) = (\lambda *\nu ) (f)$, \\
since $\lambda ^u(T_sg(u)) = \lambda ^u(g(u)) = \lambda (g)$,
particularly with $g(x)=T_uf(x)$ or  $g(u)=\nu ^x(T_uf(x))$
correspondingly, whilst $T_s\equiv {\hat v}_s$ in the cases $j=+$ or
$j=\vee $ or $j=\wedge $ with $\rho \equiv 1$, or for general
$T_uf(x)=\rho (u,x) {\hat v}_sf(x)$ in the cases of homogeneous
functionals $j=(\vee ,h)$ or $j=(\wedge ,h)$ (see \S 38 also), hence
$\nu
*\lambda , ~ \lambda *\nu \in H_j(G,K)$. Therefore, the latter formula
and Theorem 39 imply that
\par $(\nu +_i H_j(G,K))*H_j(G,K)\subset (\nu
*H_j(G,K)) +_i (H_j(G,K)*H_j(G,K))$\par $ \subset H_j(G,K)+_iH_j(G,K)\subset H_j(G,K)$ and
\par $H_j(G,K) * (\nu +_i H_j(G,K)) \subset (H_j(G,K) * \nu ) +_i (H_j(G,K)*H_j(G,K))$
\par $ \subset H_j(G,K)+_iH_j(G,K)\subset H_j(G,K)$ \\ for each $\nu \in S_j(G,K)$ and
$+_i$ corresponding to $j$,
that is $H_j(G,K)$ is the right and left ideal in $S_j(G,K)$.


\begin{thebibliography}{99}

\bibitem{bacht} Yu. A. Bachturin. {\it Basic structures of modern
alegbra} (Moscow: Nauka, 1990).

\bibitem{fell} J.M.G. Fell, R.S. Doran. {\it Representations of
$*$-algebras, locally compact groups, and Banach $*$-algebraic
bundles} (Boston: Acad. Press, 1988).

\bibitem{hew} E. Hewitt, K.A. Ross. {\it Abstract harmonic analysis}
(Berlin: Springer-Verlag, 1994).

\bibitem{kaschb} F. Kasch. {\it Moduln und Ringe} (Stuttgart:
Teubner, 1977).

\bibitem{kilknmikhmb} M. Kilp, U. Knauer, A.V. Mikhalev.
{\it Monoids, acts and categories} (Berlin: De Gruyter, 2000).

\bibitem{kunenb} K. Kunen. {\it Set theory} (Amsterdam: North-Holland
Publ. Co., 1980).

\bibitem{litmasshmz69} G.L. Litvinov, V.P. Maslov, G.B. Shpiz.
{\it Idempotent functional analysis: an algebraic approach} // Math.
Notes {\bf 69: 5-6} (2001), 696-729.

\bibitem{lujms147:3:08} S.V. Ludkovsky.
{\it Topological transformation groups of manifolds over
non-Archimedean fields, representations and quasi-invariant
measures. I} // J. Mathem. Sci. {\bf 147: 3} (2008), 6703-6846.

\bibitem{lujms150:4:08} S.V. Ludkovsky. {\it Topological transformation groups of
manifolds over non-Archimedean fields, representations and
quasi-invariant measures, II} // J. Mathem. Sci.  {\bf 150: 4}
(2008), 2123-2223.

\bibitem{mendelsb} E. Mendelson. {\it Introduction to mathematical
logic} (Princeton: D. van Nostrand Co., Inc., 1964).

\bibitem{nai} M.A. Naimark. {\it Normed rings} (Moscow: Nauka, 1968).

\bibitem{radcmcu98} T.Radul. {\it On the functor of order-preserving functionals} //
Comment. Math. Univ. of Carolinae. {\bf 39:3} (1998), 609-615.

\bibitem{rowenrb} L. Rowen. {\it Ring Theory}, Pure and
Appl. Math. Ser. {\bf 127, 128} (New York: Acad. Press, 1988).

\bibitem{scahferb} R. D. Schafer. {\it An introduction to nonassociative
algebras} (New York: Acad. Press, 1966).

\bibitem{zarichizv10} M.M. Zarichnyi. {\it Spaces of idempotent
measures} // Izvestiya Mathematics {\bf 74: 3} (2010), 481-499.

\end{thebibliography}
\end{document}